\input amstex
\documentstyle{amsppt}
\magnification=1200
\overfullrule=0pt
\NoRunningHeads
%
\def\C{\Bbb C}
\def\F{\Bbb F}

\def\R{\Bbb R}

\def\Q{\Bbb Q}
\def\Z{\Bbb Z}

%
\topmatter
\title
Three-manifolds class field theory \\ (Homology of coverings for a non-virtually Haken manifold)
\endtitle
%
\author
Alexander Reznikov
\endauthor
%
%
%
%
%
%
%
%
%
%
%
%
%
%
%
%
%
%
%
%

\toc
\head {  Introduction} \page{2}\endhead
\head {  1. Preliminary results in group cohomology} \page{4} \endhead
\head {  2. Free $\Z_p$-action in three-manifolds} \page{7} \endhead
\head {  3. Free $\Z_p \oplus \Z_p$-action in three-manifolds} \page{8} \endhead 
\head {  4. Algebraic study of linking forms} \page {12} \endhead
\head {  5. Linking forms and transfer nonvanishing theorem} \page {13} \endhead
\head {  6. The structure of anisotropic extensions, I} \page {14} \endhead
\head {  7. Shrinking} \page {16} \endhead
\head {  8. Isotropic extension} \page {17} \endhead
\head {  9. Blowing up $H_1$: the non-abelian step} \page {18} \endhead
\head {10. Class towers and the structure of the pro-p completion} \page {19} \endhead
\head {11. Nazarova-Royter theory and the structure of extensions, II} \page {22} \endhead
\head {12. Non-existence of split anisotropic constituents} \page {24} \endhead
\head {13. A spectral sequence in group cohomology} \page {26} \endhead
\head {14. Strengthened Adem inequalities and other applications} \page {27} \endhead
\head {15. Multiplication in cohomology} \page {29} \endhead
\head {16. Homology of coverings for $H_1(M)_{(2)} = \Z/2\Z \oplus \Z/4 \Z$} \page {32} \endhead
\head {17. References} \page {36} \endhead

\endtoc

\endtopmatter
\document

\newpage
$   $\newline
I am moving so slowly, \newline
A yard per five minutes \newline
How can one reach one's destination \newline
If the shoes are of small size? [O]

$$   $$\newline
This is a first in a series of papers in which we explore a deep conceptual relation between three-manifolds and number fields, a subject which may be given a name of ``arithmetic topology''. An immediate and transparent reason for such a relation is the fact that the cohomological dimension of
$Gal (\bar K|K)$
equals three for any number field $K$. \newline
This paper resulted from the author's attempt to settle the Thurston's covering conjecture - that is, that any irreducible three-manifold $M$ with infinite
$\pi_1$
has a finite covering with positive Betti number. In the process of study it became clear, however, that conjecturally non-existing non-virtually Haken manifolds show ``a strong wish to survive'', that is, the information which one can derive about them organizes itself into a harmonious picture of a non-contradictory nature. So instead of trying to show that such manifolds do not exist, we adopt a more positive approach to study the homology of their finite coverings. This is parallel to studying ideal class groups of finite extensions of a given number field. In fact, the techniques we develop will be very useful for this number-theoretic problem. \newline

The fundamental questions in three-manifold class field theory that we ask are:
\roster
\item"$\bullet$" {\it Class towers}. Suppose $M$ is a non-virtually Haken three-manifold. Fix a prime $p$. How fast the $p$-component of the first homology group
$H_1(N,\Z)$
grows for finite coverings $N$ of $M$? \newline
\item"$\bullet$" {\it Ideal class modules}. What is the structure of
$H_1(N)_{(p)}$
as a Galois module over
$G = \pi_1(M)/\pi_1(N)$
for $N$ normal? \newline
\item"$\bullet$" {\it Maximal unramified $p$-extensions}. What is a structure of pro-$p$ completion of
$\pi_1(M)$
as a pro-$p$ group?
\endroster

We also ask about the multiplication in cohomology for a virtually non-Haken three-manifold. (We give an answer in case of
$\Z/2 \Z \oplus Z/2 \Z \oplus \Z/2\Z$
a 2-torsion in
$H_1$).

A number of algebraic techniques, some well-known, some less known and some new, is needed to answer our fundamental questions. These include: \newline
\roster
\item"$\bullet$" Standard Lyndon-Serre-Hochschild spectral sequence in group cohomology \newline
\item"$\bullet$" Cohomological theory for
$\Z_p$
and
$\Z_p \oplus \Z_p$
-actions on rational homology spheres \newline
\item"$\bullet$" Linking form
$H_1 \times H_1 \rightarrow \Q/\Z$
and the classification of coverings to isotropic and non-isotropic \newline
\item"$\bullet$" Nazarova-Royter classification of f.g. modules over a cyclic group \newline
\item"$\bullet$" A new powerful spectral sequence in pro-$p$ group cohomology.
\endroster

For the first part of the paper, we will make a following ``genericity'' assumption on the three-manifold $M$, for $M$ a homology sphere:

\demo{(R) \underbar{Rich $\pi_1$}}
A three homology sphere satisfies 
$(R)$
if either the Casson invariant
$|\lambda (M)| > \sharp$
(representations of
$\pi_1(M)$
to
$SL_2(\F_5)$)
or $M$ is hyperbolic.
\enddemo

If $M$ is a not a homology sphere, we assume that
$H_1(M)_{(p)}$
has rank
$\ge 4$
for some $p$. The relation between the two assumptions is given by the following two theorems.

\proclaim {Theorem 9.1}
If $M$ satisfies the condition $(R)$ then either $M$ is virtually Haken or there exists a Zariski dense representation of
$\pi_1(M)$
in
$SL_2(\C)$.
\endproclaim

\proclaim{Theorem 9.2}
If $M$ satisfies $(R)$ then either $M$ is virtually Haken or for any prime $p$ there exists a finite normal covering of $M$ with rank
$(H_1(N,\F_p)) \ge 4$.
\endproclaim

Yet the Theorem 9.2 follows from the Theorem 9.1 and simple argument in strong approximation theory for algebraic groups as described by Lubotzky in [Lu], we describe also an alternative way using the spectral sequences in group cohomology. By the reason which will become clear later, the Theorem 9.2 represents the \underbar{non-abelian step in blowing up 
$H_1$}. The next is the abelian step.

\proclaim{Theorem 10.1 (Class towers)}
Let $p$ be a prime such that rank
$(H_1(M ,\F_p)) \ge 4$.
Suppose $M$ is not virtually Haken. Let
$M_{i+1}, \, i\ge 1$ be the maximal abelian $p$-covering of
$M_{i-1}$
(the Hilbert class covering). Let
$r_i = rank (H_1(M_i,\F_p))$.
Then
\roster
\item"(i)" $r_{i+1} \ge \frac {r^2_i -r_i} 2$. \newline
\item"(ii)" Let
$\tilde M_{i+1}$ be maximal elementary abelian $p$-covering of
$\tilde M_i \, (\tilde M_1 = M_1)$. 
The group
$H_1(\tilde M_{i+1},\Z)_{(p)}$
has exponent
$\ge p^{\tilde r_i-1}$.
In particular,
$r_i$
has superexponential growth and the class tower is infinite.
\endroster
\endproclaim

The next result answers the problem about the pro-$p$ completion of
$\pi_1(M)$.
Philosophically, it says that a non-virtually Haken 3-manifold gives rise to a ``$p$-adic 3-manifold''. On the other hand, it shows that if the Thurston conjecture is false, then there are many Poincar\'e duality pro-$p$ groups in dimension 3. (It is very hard to construct a Poincar\'e duality pro-$p$ group which is not $p$-adic analytic.)

\proclaim{Theorem 10.2}
(The structure of the pro-$p$ completion). In conditions of theorem 10.1, let
$\frak G$
be the pro-$p$ completion of
$\pi_1(M)$.
Then
$\frak G$
is a Poincar\'e duality pro-$p$ group of dimension 3.
\endproclaim

The second type of results deal with manifolds with small
$H_1(M)_{(p)}$.
For
$H_1(M)_{(p)}$
cyclic it is easy to prove that
$M \approx N/(\Z/p^n\Z)$
where $N$ is a $p$-homology sphere. For
$H_1(M)_{(2)} = \Z/2\Z \oplus \Z/2\Z$
we prove
$M \approx N/Q_{2^n}$
where $N$ is a 2-homology sphere and
$Q_{2^n}$
is the (generalized) quaternionic group. Moreover,
$n=3$
if and only if the linking form is hyperbolic. For
$H_1(M) = \Z/2\Z \oplus \Z/2 \Z \oplus \Z/2\Z$
we determine the structure of multiplication in cohomology as follows:

\proclaim{Theorem 15.4}
Let $M$ be a non-virtually Haken manifold with
$H_1(M)_{(2)} = \Z/2\Z \oplus \Z/2 \Z \oplus \Z/2\Z$.
Then there exists a basis
$x,y,z$
in
$H^1(M,\F_2)$
such that either
$xy = xz = yz =0$
or
$x^2=yz, \, y^2 =xz, \, z^2 =xy$.
\endproclaim

We then apply this analysis to determine the 2-torsion in
$\Z_2$
- and
$\Z_2 \oplus \Z_2$
-coverings of $M$ if
$H_1(M)_{(2)} \approx \Z/2\Z \oplus \Z/4\Z$.
This is the most difficult part of the paper. Its understanding requires all techniques used in the theory.

In course of our study in chapter 13, we introduce a new spectral sequence, with
$E_1$-term
$H(G,\F_p)$,
converging to
$H(K,\F_p)$
where
$1 \rightarrow K \rightarrow G \rightarrow C_p \rightarrow 1$
is an exact sequence of pro-$p$ groups. In case
$p=2$
it reduces to a long exact sequence

$$H^i(G,\F_2) \overset {res} \to \longrightarrow H^i(K,\F_2) \overset t \to \longrightarrow H^i(G,\F_2) \overset {\times s} \to \longrightarrow H^i(G,\F_2) \longrightarrow \ldots$$
which is equally important for three-manifold theory and number theory. We give some immediate applications to inequalities in group cohomology (see Adem [Ad]) and strong ring structure results about
$H^{\ast }(G,\F_2)$
in spirit of Serre's theorem [Se1], [Se2].

$$   $$
I wish to thank all people whose advice, support and insight encouraged me to device tools for entering an area of research so far from my previous experience. I wish to specially mention David Blanc, Bill Dwyer and David Kazhdan for their help. I would like to thank all participants of Bilkent Summer School in Algebraic Geometry, 1995 where I first lectured on the subject and could choose the best way of organizing the presentation. I would like also to thank the various institutions, including Universit\'e de Montreal, Max-Planck-Institut and Tibetan Library of Works and Archives, Dharamsala, where the work on the paper has been carried out. Finally, it was Lucy Katz whose unshakeable belief that the next attempt should be made after every failure just did not leave me a choice but to dig deeper. This paper is dedicated to her.

\head{1. Preliminary results is group cohomology} \endhead

In this section, we collect some miscellaneous facts from group cohomology, which will be at use in further study.

\proclaim{1.1 Proposition (Cartan)}
Let $p$ be a prime and let
$G = \undersetbrace \text{$n$} \to {\Z_p \oplus \ldots \oplus \Z_p}$.
Consider the graded ring
$\Cal A = H^*(G,\Z)$.
Then
\roster
\item"1)" $\Cal A_n$
is a vector space over
$\F_p$
for all
$n \ge 1$ \newline
\item"2)" the natural graded ring homomorphism
$\Cal A \rightarrow H^*(G,\F_p)$
is injective for all
$n \ge 1$
\endroster
\endproclaim

\demo{Proof}
Recall that
$H^*(G,\F_p) \approx S(x_1,\ldots x_n) \oplus \wedge (y_1,\ldots y_n), \, \deg x_i = 2, \, \deg y_j =1$
if $p$ is odd and
$H^*(G,\F_2) \approx S(y_1 \ldots y_n)$
if
$p=2$.
So the Poincar\'e series of
$H^*(G,\F_p)$
is
$(\frac 1 {1-t^2})^n (1+t)^n = \frac 1 {(1-t)^n}$.
Let
$K = \undersetbrace \text{$n-1$} \to {\Z_p \oplus \ldots \oplus \Z_p} < G$.
Write the Lyndon-Serre-Hochschild spectral sequence for
$H^*(G,\Z)$:

$$\searrow \quad \vdots$$
$$ \quad 0 $$
$$E_2 \qquad \qquad H^*(K,\Z_p) \Longrightarrow H^*(G,\Z)$$
$$\searrow 0$$
$$H^*(K,\Z)$$
\enddemo

We get
$P_G(t) \preceq P_K(t) + \frac {t^2} {1-t^2} \frac 1 {(1-t)^{n-1}}$,
where
$P_G(t)$
is the Poincar\'e series
$1 + \sum_{k\ge 1} \log_p(|H^k(G,\Z)|)$.
Observe that equality will imply that the spectral sequence above is degenerate at
$E_2$.
Now, the short exact sequence
$0 \longrightarrow \Z \overset {\times p} \to \longrightarrow \Z \longrightarrow \F_p \longrightarrow 0$
implies the long exact sequence

$$H^i(G,\Z) \overset {\times p} \to \longrightarrow H^i(G,\Z) \longrightarrow H^i(G,\F_p) \longrightarrow H^{i+1} (G,\Z) \longrightarrow \ldots ,$$
for
$i \ge 1$,
so
$|H^i(G,\F_p)| \le |H^i(G,\Z)||H^{i+1} (G,\Z)|$
and the equality implies that
$H^i(G,\Z)$
is a vector space over
$\F_p$.
This gives

$$\sum_{\ell \ge 1} \dim_{\F_p} (H^i(G,\F_p)) \preceq \frac {(1+t)} t (P_G(t)-1) \preceq \frac {(1+t)} t (P_K(t) -1) + \frac t {(1-t)^n}$$
Now, suppose by induction that the statement of the theorem is valid for
$(n-1)$,
then
$\frac {(1+t)} t (P_K(t)-1) = \left( \sum_{i\ge 1} \dim_{\F_p} (H^i(K,\F_p)) \right) = \frac 1 {(1-t)^{n-1}} -1$,
so
$\sum_{i\ge 1} \dim_{\F_p} (H^i(G,\F_p)) \preceq \frac 1 {(1-t)^n} -1$. \newline
Since this is actually an equality, we have equalities elsewhere above, which proves the induction step, hence the theorem.

\proclaim {1.2 Corollary}
Let
$n+2$.
Then:
\roster
\item"1)" if
$p=2$
then
$\tilde \Cal A \subset S(y_1,y_2)$
is a subring
$\F_2[y^2_1,y^2_2,y_1y_2(y_1 + y_2)]$ \newline
\item"2)" if
$p>2$
then
$\tilde \Cal A \approx (\F_p[x_1,x_2])[\epsilon]$
where
$\epsilon^2 =0$
and
$\deg \epsilon =3$.
Here
$\tilde \Cal A = \F_p \oplus_{i \ge 1} \Cal A_i$ \newline
\item"3)" $\dim_{\F_p} \Cal A^{2m} = m+1, \, \dim_{\F_p} \Cal A^{2m+1} = m$
\endroster
\endproclaim

\demo{Proof}
The s.s. above looks like

$$\matrix \format \c & \quad \c & \quad \c & \quad \c & \quad \c & \quad \c & \quad \c & \quad \c & \quad \c & \quad \c \\
 &  & &  \ldots \\
\text{   }\\
 &  \Z_p & \Z_p & \Z_p & \Z_p & \Z_p & \Z_p & \Z_p & 0 & \vdots \\
\text{   }\\
E_2 & 0 & 0 & 0 & 0 & 0 & 0 & 0 & 0 & \vdots \\
\text{   }\\
& \Z_p & \Z_p & \Z_p & \Z_p & \Z_p & \Z_p & \Z_p & \Z_p &\vdots \\
\text{   } \\
& 0 & 0 & 0 & 0 & 0 & 0 & 0 & 0 & \vdots \\
\text{   }\\
& \Z & 0 & \Z & 0 & \Z & 0 & \Z & 0
\endmatrix$$
so
$\Cal A^3 \equiv \F_p$.
In 1) the only element in
$\F^{(3)}_2(y_1,y_2)$,
invraint under the natural action of
$GL_2(\F_2)$
is
$y_1y_2(y_1 + y_2)$.
2), 3) an obvious from the spectral sequence.
\enddemo

\proclaim{1.3 Proposition}
Let $W$ be a
$\Z_p$
-module, finitely generated as an abelian group. Then:
\roster
\item"1)" $\oplus \Sb i \, \text{odd} \\ i \ge 1 \endSb H^i(\Z_p,W)$
is a free
$\F_p[t]$
module of rank
$\dim_{\F_p}H^1(\Z_p,W)$ \newline
\item"2)" $\oplus \Sb i \, \text{even} \\ i \ge 2 \endSb H^i(\Z_pW)$
is a free
$\F_p[t]$
module of rank
$\dim_{\F_p} H^2(\Z_p,W)$
\endroster
where the module structure comes from the natural coupling
$H^*(Z_p,W) \otimes H^*(\Z_p,\Z) \rightarrow H^*(\Z_p,W)$.
\endproclaim

\demo{Proof}
See [CE].
\enddemo

\proclaim{1.4 Proposition}
Let $G$ be a group acting freely and discontinuously in the connected sphere $X$. Let
$Y = X/G$.
Consider the 1.1. of the covering
$X \rightarrow Y$:

$$\qquad \searrow \qquad H^*(G,H^2(Y,\Z)$$

$$E_2 \qquad \qquad H^*(G,H^1(Y,\Z))$$

$$\searrow \qquad H^*(G,\Z)$$
Then at any
$E_r$
all rows are graded
$H^*(G,\Z)$
modules and all differentials are module homomorphism.
\endproclaim

\demo{Proof}
This follows from the multiplicative structure of the equivalent s.s. of the Borel's fibration

$$\CD
X @>>> X_G \\
@. @VVV \\
@. BG \endCD$$
\enddemo

\demo{1.5}
Let $W$ be a module over
$\Z_p \oplus \Z_p$.
Consider the
$E_2$
-term of the Lyndon-Serre-Hochschild spectral sequence
$H^i(\Z_p,H^i(\Z_p,W)) \Rightarrow H^{i+1} (\Z_p \oplus \Z_p,W)$.
There is a natural action of the ring without unit
$(X,Y) \subset \Z_p[X,Y]$
in all
$E_r$
with the properties:
\roster
\item"1)" multiplication by
$X : E^{p,q}_2 \rightarrow E^{p+2,q}_2$
is an isomorphism for
$p > 0$ \newline
\item"2)" multiplication by
$Y : E^{p,q}_2 \rightarrow E^{p,q+2}_2$
is an isomorphism for
$q > 0$ \newline
\item"3)" $d_r$
is an endomorphism of 
$E_r$
as a graded
$(X,Y)$
-module, \newline
\item"4)" the action of
$(X,Y)$
in
$E_{\infty}$
agrees with the module structure of
$H^*(\Z_p \oplus \Z_p,W)$
our
$H^*(\Z_p \oplus \Z_p\Z)$.
\endroster
\enddemo

\demo{Proof}
The proof is immediate from the double complex, which is a resolution of $\Z$ over the sum of two cyclic groups.
\enddemo

\head{2. Free $\Z_p$-actions in three-manifolds} \endhead

This and the next section are devoted to the cohomological study of elementary abelian group actions on rational homology three-spheres.

\demo{2.1}
Suppose
$\Z_p$
acts freely in a closed oriented three-manifold $N$, which is a rational homology sphere, that is,
$H_1(N,\Z)$
is torsion. Then
$W = H_1(N,\Z)$
is a
$\Z_p$
-module. We wish to understand the cohomology
$H^i(\Z_p,W)$.
Put
$M = N/\Z_p$
\enddemo

\proclaim{Proposition 2.1}
\roster
\item"1)" either
$\dim_{\F_p}H^1(\Z_p,W)$
and
$\dim_{\F_p} H^2(\Z_p,W)$
are both 1,
\endroster
or they are both zero.
\roster
\item"2)" there is a natural exact sequence
$0 \rightarrow \F_p \rightarrow H^2(M,\Z) \rightarrow H^0(\Z_p,W) \rightarrow 0$.
\endroster
\endproclaim

\demo{Proof}
Write the cohomological s.s. for
$N \rightarrow M$,
observing that
$H^2(N,\Z) \approx W$
by Poincar\'e duality

$$\matrix \format \c \quad & \c \quad & \c \quad & \c \quad & \c \quad & \c \quad & \c \quad & \c \quad & \c \quad & \c \quad &\c \\
& \Z & 0 & \Z_p & 0 & \Z_p & 0 & \Z_p \\
\text{   }\\
E_2 & H^0(\Z_p,W) & H^1 &H^2 & H^3 &H^4 &H^5 &H^6 &\Rightarrow & H^*(M,\Z) \\
\text{   }\\
& 0 & 0 & 0 & 0 & 0 & 0 & 0\\
\text{   }\\
& \Z & 0 & \Z_p & 0 & \Z_p & 0 & \Z_p
\endmatrix$$
Since
$H^k (M,\Z) =0$
for
$k \ge 4$
all rows should be killed. That means that either the
$d_2$
kills all odd
$H^i(\Z_p,W)$
and
$d_3$
kills all even
$H^i(\Z_p,W)$
or
$d_4$
kills the first and the fourth row by Proposition 1.3 and Proposition 1.4. That proves 1), and 2) is obvious.
\enddemo

\head{3. Free $\Z_p \oplus \Z_p$-actions in three-manifolds} \endhead

\demo{3.1}
Suppose now that a rational homology sphere
$N^3$
is acted upon freely by
$\Z_p \oplus \Z_p$.
Then
$W = H^2(N,\Z)$
is a
$\Z_p \oplus \Z_p$
module, and we are interested to understand
$H^*(\Z_p \oplus \Z_p,W)$.
Write the group extension s.s. [CE] using the information of Proposition 2.1:

$$ \matrix \format \c \quad & \c \quad & \c \quad & \c \quad & \c \quad & \c \quad & \c \quad & \c \quad & \c \quad  & \c\\
&\Z_p & \Z_p & \Z_p \\
\text{   }\\
& \Z_p & \Z_p & \Z_p \\
\text{   }\\
E_2 & \Z_p & \Z_p & \Z_p & \Rightarrow & H^*(\Z_p \oplus \Z_p,W)
\text{   } \endmatrix$$
$\qquad \qquad \qquad \qquad \qquad \quad H^*(\Z_p,H^0(\Z_p,W))$\newline
$   $\newline
or

$$ \matrix \format \c \quad & \c \quad & \c \quad & \c \quad & \c \quad & \c \quad & \c \quad & \c \quad & \c \quad  & \c\\
& 0 & 0 & 0 & 0 \\
\text{  }\\
E_2 & 0 & 0 & 0 & 0 &\Rightarrow &H^*(\Z_p \oplus \Z_p,W) \\
\text{   } \endmatrix$$
$\qquad \qquad \qquad \qquad \qquad \quad H^*(\Z_p,H^0(\Z_p,W)$ \newline
$   $\newline
We will call this alternative case A and case B. In any case, we wish to understand the first row, namely,
$H^*(\Z_p,(H^0(\Z_p,W))$.
By Proposition 2.1 2), we have a short sequence of
$\Z_p$
-modules
$0 \rightarrow \F_p \rightarrow H^2(M,\Z) \rightarrow H^0(\Z_p,W) \rightarrow 0$.
Now, the first factor of
$\Z_p \oplus \Z_p$
acts freely on
$M = N/$
(action of the second factor), and $M$ is a rational homology sphere, because $N$ is (recall that
$H^*(M,\Q) \approx H^*_{inv} (N,\Q))$.
So by proposition 2.1 1), either
$H^i(\Z_p,H^2(M,\Z)) = \Z_p$
or 0 for all
$i \ge 1$.
Now, the long exact sequence

$$\ldots \rightarrow H^i(\Z_p,\F_p) \rightarrow H^i(\Z_p,H^2(M,\Z)) \rightarrow H^i(\Z_p,H^0(\Z_p,W)) \rightarrow H^{i+1} (\Z_p,\F_p) \rightarrow \ldots$$
reduces to either
$H^i(\Z_p,H^0(\Z_p,W)) \approx H^{i+1} (\Z_p,\F_p) \approx \F_p$,
or to
$\ldots \F_p \rightarrow \F_p \rightarrow H^i(\Z_p,H^0 (\Z_p,W)) \rightarrow \F_p \rightarrow \F_p \rightarrow H^{i+1} (\Z_p,H^0(\Z_p,W)) \rightarrow \ldots$ \newline
Now, in the latter case the map
$H^*(\Z_p,\F_p) \rightarrow H^*(\Z_p,H^2(M,\Z))$
is a
$\F_p[t]$
-module homomorphism
$(\deg t = 2)$,
so it is either zero or an isomorphism for all $i$ of the same parity. This implies immediality that all
$H^i(\Z_p,H^0(\Z_p,W))$
are of the same dimension 0,1, or 2 for
$i \ge 1$. \newline
So the
$E_2$
term for
$H^*(\Z_p \oplus \Z_p,W)$
looks like:

$$\text{Case A}$$

$$\matrix \format \c \quad & \c \quad & \c \quad & \c \quad & \c \quad & \c \quad & \c \quad & \c \quad & \c \quad & \c\\
& \Z_p & \Z_p & \Z_p & \Z_p \\
\text{   }\\
E_2 & \Z_p & \Z_p &\Z_p &\Z_P & 0 \le m \le 2, &\Rightarrow &H^*(\Z_p \oplus \Z_p,W) \\
\text{   } \\
& * &\Z^m_p & \Z^m_p & \Z^m_p & \ldots ,
\endmatrix$$

$$\text{Case B}$$

$$\matrix \format \c \quad & \c \quad & \c \quad & \c \quad & \c \quad & \c \quad & \c \quad & \c \quad & \c \quad & \c\\
& & & 0 \\
\text{   }\\
E_2 &* &\Z^m_p &\Z^m_p &\Z^m_p &\ldots & 0 \le m \le 2.
\endmatrix$$
In particular,
$\dim_{\F_p} H^{\ell}(\Z_p \oplus \Z_p,W)) \le m+\ell$
(case A) and
$\le m$
(case B).
\enddemo

\demo{3.2}
Let
$Q =N/\Z_p \oplus \Z_p$
and write the s.s. for the covering
$N \rightarrow Q$:

$$H^*(\Z_p \oplus \Z_p,\Z) \approx \Cal A$$

$$E_2 \quad \searrow \quad H^*(\Z_p \oplus \Z_p,W) \quad \rightarrow H^*(Q,\Z)$$

$$0$$

$$H^*(\Z_p \oplus \Z_p,\Z) \approx \Cal A$$
Since
$H^i(M,\Z) =0$
for
$i \ge 4$,
all rows should be killed by differentials in high dimensions. There are exactly three nontrivial differentials which are 
$\Cal A$
-homomorphisms:
$d_2 : \Cal A \rightarrow H^*(\Z_p \oplus \Z_p,W), \, d_3 : H^*(\Z_p \oplus \Z_p,W)/Im \, d_2 \rightarrow \Cal A$
and
$d_4 : \text{Ker}\, d_2 \rightarrow \Cal A/Im \, d_3$. \newline
We first observe:
\enddemo

\proclaim{Proposition 3.1}
$d_3$
is nontrivial in high dimension.
\endproclaim

\demo{Proof}
Suppose the opposite. Then
$d_4 : \text{Ker} \, d_2 \rightarrow \Cal A_{n+4}$
should be an isomorphism, which is impossible by Corollary 1.2.3).
\enddemo

\proclaim{Proposition 3.2}
$d_4 =0$
in positive dimensions.
\endproclaim

\demo{Proof}
Suppose the opposite. Assume first
$p=2$.
Then
$0 \ne \text{Ker} \, d_2 \approx \Cal A/Im \, d_3$.
$\text{Ker} \, d_2$
is an ideal in
$\Cal A$,
so it contains a principal ideal and since
$\Cal A$
is a domain by Proposition 1.1, 1.2.1), the dimension of
$(\text{Ker} \, d_2)_n$
grows as $n$. On the other hand,
$Im \, d_3$
is an ideal in
$\Cal A$.
We may look at
$\Cal A$
as a graded ring of an irreducible projective curve over
$\F_2$.
Then
$\Cal A/Im \, d_3$
is a graded ring of a null-dimensional scheme, so
$\dim_{\F_2}(\Cal A/Im \, d_3)$
stays bounded, a contradiction. \newline
Now, let $p$ be odd. Recall
$A = \F_p[x,y][\epsilon], \, \epsilon^2 =0$.
Any homogeneous ideal $I$ contains
$(f) \cdot \epsilon$,
where
$0 \ne f \in \F_p[x,y]$,
and so
$\dim I_{2n-1}$
grows as
$2n$. So
$\dim (\text{Ker} \, d_2)_{2n-1} \approx \dim (\Cal A/Im \, d_3)_{2n+3}$
grows as
$2n$.
On the other hand again,
$Im \, d_3$
is nontrivial in high dimension, so
$\dim (\Cal A/Im \, d_3)_{2n+3}$
stays bounded by the same argument as above. So
$d_4 =0$. \newline
Observe that
$E^{\infty}$
term of the spectral sequence above should be

$$\matrix \c \quad \format \c \quad & \c \quad & \c \quad & \c \quad & \c \quad & \c \quad & \c \quad & \c \quad & \c \quad & \c\\
& p^2\Z & 0 & 0 & 0 \\
\text{   } \\
E^{\infty} & 0 & 0 & 0 & 0\\
\text{   }\\
& 0 & 0 & 0 & 0\\
\text{   }\\
& \Z & 0 &\Z_p\oplus \Z_p
\endmatrix$$
since
$H^3(M,\Z) \approx \Z$
does not have torsion. \newline
That implies the following.
\enddemo

\proclaim{3.3 Corollary}
For
$n \ge 3$
we have an exact sequence
$0 \rightarrow \Cal A_{n-2} \overset {d_2} \to \longrightarrow H^n (\Z_p \oplus \Z_p,W)\overset {d_3} \to \longrightarrow \Cal A_{n+3} \rightarrow 0$. \newline
In particular,

$$\align a_{2n} = &\dim_{\F_p} H^{2n} (\Z_p \oplus \Z_p,W) = n + (n+1) = 2n+1 \\
a_{2n-1} = &\dim_{\F_p} H^{2n-1} (\Z_p \oplus \Z_p,W) = (n-2) + (n+2) = 2n, \, n\ge 1 \\
&\dim_{\F_p} H^1(\Z_p \oplus \Z_p,W) \le 3, \, \dim H^2(\Z_p \oplus, \Z_p, W) \le 4. \endalign$$
\endproclaim

\demo{4.5}
Now recall the s.s. for
$H^*(\Z_p \otimes \Z_p,W)$
from section 2. We see immedaitely that case B is imossible because of the inequality
$\dim (H^2(\Z_p \oplus \Z_p,W)) \le m$. 
So we get the s.s.

$$\matrix \format \c \quad & \c \quad & \c \quad & \c \quad & \c \quad & \c \quad & \c \quad & \c \quad & \c \quad & \c\\
& \Z_p & \Z_p & \Z_p &\Z_p \\
\text{   }\\
E_2 & \Z_p & \Z_p &\Z_p &\Z_p &\Rightarrow &H^*(\Z_p \oplus \Z_p,W) \\
\text{   }\\
& * &\Z^m_p &\Z^m_p &\Z^m_p &\ldots
\endmatrix$$
and
$a_{\ell} = \dim_{\F_2} H^{\ell}(\Z_p \oplus \Z_p,W) \le \ell +m$.
\enddemo

In particular, since for big
$n, \, a_{2n-1} =2n$,
we have
$m \ge 1$,
so
$m =1$
or
$m=2$.
Assume
$m=2$.
If
$a_2 =3$
there should be a nontrivial differential, killing something from the second diagonal
$(\Z_p,\Z_p,\Z^2_p)$.
We have the following choices

$$\matrix \format \c \quad & \c \quad & \c \quad & \c \quad & \c \quad & \c \quad & \c \quad & \c \quad & \c \quad & \c\\
\Z_p & \Z_p & \Z_p \\ 
\text{   }\\
\Z_p & \Z_p & \Z_p \\
\text{   }\\
& d_2 \ne 0 \\
\text{   }\\
* &\Z^2_p &\Z^2_p
\endmatrix \tag case I$$

$$\matrix \format \c \quad & \c \quad & \c \quad & \c \quad & \c \quad & \c \quad & \c \quad & \c \quad & \c \quad & \c\\
\Z_p & \Z_p & \Z_p \\
\text{   }\\
& d_2 \ne 0 \\
\text{   }\\
\Z_p & \Z_p & \Z_p \\
\text{   }\\
* &\Z^2_p &\Z^2_p
\endmatrix \tag case II$$

$$\matrix \format \c \quad & \c \quad & \c \quad & \c \quad & \c \quad & \c \quad & \c \quad & \c \quad & \c \quad & \c\\
\Z_p & \Z_p & \Z_p \\ 
\text{   }\\
\Z_p & \Z_p & \Z_p \\
\text{   }\\
*\Z^2_p &\Z^2_p &\Z^2_p & \text{either} \quad d_2 \text{or} \quad d_3 \ne 0
\endmatrix \tag case III$$
If the case I is not realized, then
$a_1 =\dim H^1(\Z_p \oplus \Z_p,W) =3$,
but that means that in the spectral sequence of page 6,
$d_3 : H^1 (\Z_p \oplus \Z_p,W) \rightarrow A_4$
is an isomorphism, so
$d_4 : \text{Ker} \, d_3  \rightarrow A_4$
is zero. Because of the form of
$E^{\infty}$
that means that \newline
$|A_0 : \text{Ker} \, d_2| = p^2$,
so
$\dim Im \, d_2 =2$,
hence
$\dim H^2(\Z_p \oplus \Z_p,W) = 2+ \dim A_5 =4$,
a contradiction. Se we have case I. Moreover, since
$a_1 =2$
and
$a_2 =3, \, d_4 : p\Z \rightarrow A_4$
is nontrivial and
$\dim (Im \, d_2) =1$.
If
$a_2 =4$,
then there may be no differential, touching the second diagonal, in particular
$a_1 =3$
and
$|Im \, d_2| = p^2$
as above. 

\head{4. Algebraic study of linking forms} \endhead

Let $W$ be a finite abelian group. A linking form is a symmetric bilinear form

$$(\cdot , \cdot) : W \otimes_{\Z} W \rightarrow \Q/\Z$$
for which an induced map

$$W \rightarrow Hom (W,\Q/\Z) = \hat W$$
is an isomorphism.

\proclaim{Lemma (4.1)}
If
$W = \oplus_p W_{(p)}$
be the canonical decomposition as a direct sum of $p$-groups, then
$W_{(p)}$
is orthogonal to
$W_{(q)}$
for
$p \ne q$.
\endproclaim

The proof is obvious. From now on we assume $W$ to be a $p$-group.

\proclaim{Lemma (4.2)}
If
$x \in W,\, ord \, (x) = p^k$,
then
$\max_y \, ord \, (x,y) = p^k$.
\endproclaim

The proof is again obvious. \newline
Let
$V = p^{\ell} W$
and define for
$u,v \in V$
\TagsOnRight

$$(u,v)_V = p^{\ell} (x,y) = (u,y) = (x,v) \tag*$$
where
$p^{\ell} x =u$
and
$p^{\ell}y =u$.

\proclaim{Lemma (4.3)}
The formula (*) defines a linking form in $V$.
\endproclaim

\demo{Proof} 
If
$\bar x \ne x$
with
$p^{\ell} \bar x =u$,
then
$p^{\ell} (\bar x,y) -p^{\ell}(x,y) = (u,y)-(u,y) =0$,
so
$(u,v)_V$
is well-defined. On the other hand, if
$(u,v)_V =0$
for all $v$, then
$(u,v) =0$
for all $y$, so
$u=0$.
\hfill Q.E.D.
\enddemo

Let
$W = \Z/p^{k_1}\Z \oplus \Z/p^{k_2} \Z \oplus \ldots \oplus \Z/p^{k_s}\Z$
with
$k = k_1 \ge k_2 > \ldots \ge k_s$.
Let
$U = W/\text{Ker}$
(multiplication by
$p^{k-1}$).
Then we have an isomorphism of
$\F_p$
-vector spaces

$$U \overset{p^{k-1}} \to \longrightarrow p^{k-1} W,$$
and
$(\cdot , \cdot)_{p^{k-1}W}$
induces a nondegenerate 
$\F_p$
-valued scalar product in $U$. If $p$ is odd there exists an element
$\bar x \in U$
with
$(\bar x,\bar x)_U \ne 0$.
Let $x$ be an element of $W$ of maximal order
$p^k$,
which projects to
$\bar x$,
then
$p^{k-1} (x,x) \ne 0$,
so
$(x,x)$
has order
$p^k$.
Let
$W_x$
be a cyclic group, generated by
$x_1$
then we have just proved that
$(\cdot , \cdot)|_{W_x}$
is nondegenerate, so $W$ splits as
$W_x \oplus (W \ominus W_x)$.
This result in the following proposition

\proclaim{Proposition (4.4)}
For $p$ odd, there exists an \underbar{orthogonal} decomposition

$$W = \Z/p^{k_1} \Z \oplus \ldots \oplus \Z/p^{k_s} \Z$$
\endproclaim

In the following lemma we assume $p$ is odd.

\proclaim{Lemma (4.5)}
Suppose
\roster
\item"1)" $W$ has a linking form \newline
\item"2)" $W$ admits an orthogonal action of
$C_p$ 

$   $\newline
Then for $\zeta$ a generator of
$C_p$-action,
$\dim_{\F_p} Im (1-\zeta)$
is even.
\endroster
\endproclaim

\demo{Proof}
Consider a pairing on
$Im (1-\zeta)$,
defined by
$\langle (1-\zeta)x,(1-\zeta)y \rangle = (x,\zeta^{\frac {p-1} 2} (1-\zeta)y)$.
It is immediately seen to be nondegenerate. Moreover, since
$\zeta$
acts orthogonally, we have

$$\langle (1-\zeta)y,(1-\zeta)x \rangle = (y,\zeta^{\frac {p-1} 2}(1-\zeta)x) = (\zeta^{\frac {p+1} 2}(1-\zeta^{-1}) y,x) = -\langle (1-\zeta)x, (1-\zeta)y \rangle,$$
so
$\langle \cdot \, , \cdot \rangle$
is a symplectic structure. Arguing like in 1.4, we see that
$W \approx \sum_i W_i \oplus \hat W_i$
so rank $W$ is even and
$\dim_{\F_p}W$
is even, too.
\enddemo

\head{5. Linking forms and transfer non-vanishing theorem} \endhead

\demo{5.1 \underbar{Fundamentals}} 
Let
$M^{2m-1}$
be a closed oriented manifold. The linking form (see, for example, [vD])

$$(\quad , \quad ) : H^{tors}_{m-1} (M) \otimes_{\Z} H^{tors}_{m-1} (M) \rightarrow \Q/\Z \tag 1$$
is defined as follows: by universal coefficients formula one has
$H^m_{tors} (M) \approx Ext^1 (H^{tors}_{m-1}(M),\Z) \approx Hom (H^{tors}_{m-1} (M),\Q/\Z)$.
On the other hand,
$H^m_{tors}(M) \approx H^{tors}_{m-1}(M)$
by Poincar\'e duality, so one gets an isomorphism

$$H^{tors}_{m-1} (M) \approx Hom (H^{tors}_{m-1} (M), \Q/\Z),$$
which means that there is a nondegenerate form (1), which is easily seen to be symmetric. \newline
Here is a more geometric description. Let
$x,y \in H^{tors}_{m-1} (M)$,
so that
$N \cdot x =0$.
Realize $x$ by a smooth chain
$\tilde x$
and find a smooth chain
$\tilde z$
such that
$\partial \tilde z = N \tilde x$.
Realize $y$ by a smooth chain
$\tilde y$,
disjoint from $x$ and intersecting $z$ transversally. Let
$\sharp(\tilde z \cap \tilde y)$
be a number of intersection points, counted with sign. Then

$$(x,y) = \frac 1 N \sharp (\tilde z \cap \tilde y) \quad (mod \, \Z)$$
\enddemo

\demo {5.2 Reciprocity formula}
Let
$N \overset {\pi} \to \longrightarrow M$
be a finite covering of
$2m-1$
-dimensional manifolds. Let
$t : H_i(M) \rightarrow H_i(N)$
be the transfer map. The reciprocity formula reads: for any
$x \in H^{tors}_{m-1} (N), \, y \in H^{tors}_{m-1} (M)$,

$$(\pi_*x,y)_M = (x,ty)_N,$$
or, in other words,
$t =- (\pi_*)^*$,
as linear maps between abelian groups with nondegenerate
$\Q/\Z$
-valued scalar product.
\enddemo

\demo{Proof}
Let
$\partial \tilde z = N \cdot \tilde y$,
then
$\partial(t\tilde z) = N \cdot t \tilde y$.
Realize $x$ by a smooth chain disjoint from
$t \tilde y$
and transversal to
$t \tilde z$.
Then
$\pi_*\tilde z$
is a smooth chain, disjoint from
$\tilde y$,
and transversal to
$\tilde z$.
\enddemo

Next,
$\sharp(\pi_* \tilde z,\tilde y) = \sharp (\tilde z,t\tilde y)$,
which proves the formula. \newline
One notices that the reciprocity formula follows from the formula
$t = PD \circ \pi_* \circ PD$
for the Gysin homomorphism, c.f. [Ka].

\proclaim{Transfer non-vanishing Theorem (5.3)}
Suppose
$H_i(N) =0$
for
$0 < i < m-1,\, H_{m-1}(N)$
is torsion,
$\pi : N \rightarrow M$
is a normal covering with the Galois group $G$ and suppose
$H_{m-1}(G) =0$.
Then
$t : H_{m-1}(M) \rightarrow H_{m-1}(N)$
is injective.
\endproclaim

\demo{Proof}
Write the homological spectral sequence of the covering
$N \rightarrow M$:

$$\align
(H_{m-1}&(N)_G \\
0 \quad\\
\Z \quad &H_1(G,\Z) \ldots H_{m-1} (G,\Z) \quad H_m(G,\Z) \ldots \\
& \qquad \qquad \qquad \parallel \\
& \qquad \qquad \qquad \, 0
\endalign$$
We see that
$H_{m-1} (M) \approx (H_{m-1}(N))_G$,
in particular,
$H_{m-1}(N) \overset {\pi_*} \to \longrightarrow H_{m-1} (M)$
is onto, hence
$t = (\pi_*)^*$
is injective, since
$( \qquad )_M$
and
$( \qquad )_N$
are nondegenerate.
\enddemo

\proclaim{Corollary (5.4)}
Let
$N \rightarrow M$
be a normal covering of three-manifolds with the Galois group $G$. If
$b_1(M) = b_1(N) =0$
and
$G/[G,G] =1$,
then
$t : H_1(M) \rightarrow H_1(N)$
is injective. In particular, this is true for
$G = SL(2,\F_p)$
for
$p \ge 5$.
\endproclaim

\head{6. The structure of anisotropic extension, I} \endhead

\demo{6.1}
Let $M$ be a closed oriented three-manifold with
$H_1(M)_{(p)} = \Z/p^{k_1} \Z \oplus \ldots \oplus \Z/p\Z$,
where the decomposition is orthogonal with respect to the linking form. Fix a generator $z$ of the last component. Conisder the homomorphism

$$(z,\cdot) : H_1(M) \rightarrow \Z_p \tag*$$
and denote by $N$ the cyclic covering of $M$ with respect to this homomorphism. We assume
$b_1(M) = b_1(N) =0$.
In this chapter we study the structure of
$(H_1(N))_{(p)}$.
Let
$e_1, \ldots, e_{s-1}$
be the generators of all components except the last. Then
$(e_i,e_j) =0$
and
$(e_i,e_i)$
has the order
$p^{k_i}$
in
$\Q/\Z$.
Since all
$e_i$
lie in the kernel of (*) the transfer
$te_i$
splits as
$v_i + \zeta v_i + \ldots + \zeta^{p-1}v_i$,
where
$\pi_*(v_i) =e_i$
and $\zeta$ is the deck transformation corresponding to (*). Indeed, we may find a loop in
$\pi_1(M)$,
representing
$e_i$,
and lying in the kernel of the composition
$\pi_1(M) \rightarrow \pi_1(M) \rightarrow \Z_p$,
so that its preimage in $N$ splits to $p$ connected components. We denote
$W = H_1(N)_{(p)}, \, V = H_1(M)_{(p)}$
and
$W_i$
the subgroup of $W$ generated by
$v_i,\zeta v_i \ldots \zeta^{p-1}v_i$.
We start with the following (obvious) remark.
\enddemo

\proclaim{Lemma (6.1)}
$W_i$
is a cyclic
$\Z [C_p]$
-module, i.e. the quotient of
$\Z[C_p]$.
\endproclaim

\proclaim{Lemma (6.2)}
The element
$p_i = v_i + \zeta v_i + \ldots + \zeta^{p-1}v_i = te_i$
is of order
$p^{k_i}$
in
$W_i$.
\endproclaim

\demo{Proof}
On one hand, the order of
$p_i = te_i$
is no more that the order of
$e_i$,
that is,
$p^{k_i}$.
On the other, by the reciprocity formula of 5.2,

$$(p_i,v_i)_N = (e_i,e_i)_M,$$
and the RHS is of order
$p^{k_i}$,
so
$ord (p_i) \ge p^{k_i}$,
as claimed.
\enddemo

\proclaim{Lemma (6.3)}
The restriction of
$t : H_1(M) \rightarrow  H_1(N)$
on
$\oplus^{s-1}_{i=1} \Z/p^{k_i} \Z$
is injective. Moreover
$t (\oplus^{s-1}_{i=1} \Z/p^{k_i}\Z) =(H_1(N))^{\Z_p}$.
\endproclaim

\demo{Proof}
The first statement follows from
$(p_i,v_j)_N = (e_i,e_j)_M$
th orthogonality of the decomposition of
$H_1(M)$
and nondegeneracy of
$(\cdot , \cdot)_M$.
To prove the second, write the homology spectral sequence of the covering
$N \rightarrow M$,

$$\align
&\vdots \\
\text{   }\\
&(H,(N)_{\Z_p} \qquad d_2 \\
\text{   }
& \Z \qquad \Z_p \qquad 0 \endalign$$
from which we deduce the exact sequence

$$0 \rightarrow (H_1(N))_{\Z_p} \rightarrow H_1(M) \rightarrow \Z_p \rightarrow 0$$
so
$\dim_{\F_p} H_1(N)_{\Z_p} = k_1 + \ldots + k_{s-1}$,
hence
$\dim_{\F_p} (H_1(N))^{\Z_p} = k_1 + \ldots + k_{s-1}$,
and so the injectivity of $t$ implies the surjectivity.
\enddemo

\proclaim{Corollary (6.3)}
Let
$N \rightarrow M$
be a cyclic covering, corresponding to the homomorphism
$(\quad , z): H_1(M) \rightarrow \Z_p$
with
$(z,z) \ne 0$.
Then
$H^i(\Z_p,H_1(N)) =0$
for
$i \ge 1$.
\endproclaim

\demo{Proof}
Since
$\oplus^{s-1}_{i=1} \Z/p^{k_i}\Z = \pi_*(H_1(N))$
we have
$(H_1(N))^{\Z_p} = t(\oplus^{s-1}_{i=1} \Z/p^{k_i} \Z) = t \circ \pi_*(H_1(N)) = Im (1+ \zeta + \ldots + \zeta^{p-1})$
so that
$H^1(\Z_p,H_1(N)) =0$.
But then
$H^i(\Z_p,H_1(N)) =0$
for all
$i \ge 1$
by proposition 2.1.
\enddemo

\proclaim{Lemma (6.4)}
The restriction of
$(\cdot ,\cdot)_N$
on
$W_1 + \ldots + W_{s-1}$
is nondegenerate.
\endproclaim

\demo{Proof}
Let
$y \in W_1 + \ldots + W_{s-1}$.
There exists a minimal $N$ such that
$(1-\zeta)^{N+1} y=0$,
so
$0 \ne x = (1-\zeta)^Ny$
and
$x \in (H_1(N))^{\Z_p}$.
By the previous lemma,
$x = \alpha_1p_1 + \ldots + \alpha_sp_s$
such that at least for one $i$,
$\alpha_i \not\in p^{k_i} \Z$
by lemma 3.2, Now

$$(x,v_i)_N = (\alpha_1p_1 + \ldots + \alpha_sp_s,v_i)_N = \alpha_i(p_i,v_i)_N = \alpha_i(e_i,e_i)_M \ne 0.$$
But

$$(x,v_i)_N = ((1-\zeta)^Ny,v_i)_N = (y,(1-\zeta^{-1})^Nv_i)_N,$$
and
$(1-\zeta^{-1})^N v_i \in W_i$,
so
$(\cdot , \cdot)_N$
is nondegenerate on
$W_1+ \ldots + W_{s-1}$.
\enddemo

\demo{Corollary (6.5)}
$H_1(N) = W_1 + \ldots + W_{s-1}$
\enddemo

\demo{Proof}
By the previous lemma,
$H_1(N)$
splits as
$(W_1+ \ldots + W_S) \oplus (W_1 + \ldots + W_S)^{\perp}$.
Suppose the latter space is nontrivial. Since any action of
$C_p$
in an abelian $p$-group has fixed points, we would have
$H_1(N)^{\Z_p}$
strictly contains
$t(H_1(M))$,
which contradicts lemma 3.3.
\enddemo

\proclaim{6.6 Theorem (Structure theorem for anisotropic extension)}
The group
$H_1(N)$
is a sum of cyclic modules:
$H_1(N) = W_1 \oplus \ldots \oplus W_{s-1}$
with
$W_i \approx \Lambda/\alpha_i$
where
$\Lambda = \Z[\zeta]/(\zeta^p-1) = \Z[C_p]$
and
$\alpha_i$
ideals in
$\Lambda$.
Moreover
$(\Lambda/\alpha_i) \underset {\Lambda} \to {\otimes} \Z \approx \Z/p^{k_i} \Z$,
and
$Ext^1_{\Lambda} (\Z,W) =0$.
\endproclaim

\demo{Proof}
We need only to check the last statement. It follows from the lemma 6.3 that
$H^{\text{even} \, > 0} (\Z_p,H_1(N)) =0$,
so by 2.1, also
$H^{odd}(\Z_p,H_1(N)) =0$,
hence
$H^{>0} (\Z_p,W_i) =0$.
\enddemo

\head{7. Shrinking} \endhead

Shrinking is a process leading to a
$\Z/p\Z$
split component in
$H_1$
of a specially chosen covering, as described below.

\demo{7.1}
Assume that in the decomposition
$H_1(M) = \Z/p^{k_1} \Z \oplus \ldots \oplus \Z/p^{k_s} \Z, \, k_s \ge 2$
and consider a map
$H_1(M) \overset {\epsilon} \to \longrightarrow \Z/p\Z$
sending
$e_s$
to a generator and
$e_i$
to zero,
$i < s$.
Let $N$ be the corresponding covering. Let
$v_i, i<s$
be defined as above and let
$v_s = t(e_s)$.
Again denote by
$W_i$
the subgroup of
$H_1(N)$
generated by
$v_i$
as
$\Lambda$
-module. Since
$v_s$
is invariant, in fact
$W_s$
is cyclic of order
$\le p^{k_s}$.
We claim that in fact that
$(ord \, v_s) = p^{k_s-1}$.
Indeed,
$(v_s,v_s)_N = (e_s,\pi_*v_s)_M = p(e_s,e_s)$,
so
$ord (v_s,v_s)_N = p^{k_s-1}$
and so
$ord \, v_s \ge p^{k_s-1}$.
On the other hand, for any
$z \in H_1(N),\, (v_s,z)_N = (e_s,\pi_* (z))$
and since
$\pi_*(z) \in \text{Ker} \epsilon$
we see that
$ord (e_s,\pi_*(z)) \le p^{k_s-1}$,
so
$ord \, v_s \le p^{k_s-1}$.
Hence
$W_1 \approx \Z/p^{k_s-1}$.
Now we claim that one has an exact sequence

$$0 \rightarrow \Z_p \overset{\psi} \to \longrightarrow H_1(M) \overset t \to \longrightarrow H_1(N)^{\Z_p} \rightarrow 0,$$
where
$Im \, \psi = p^{k_s-1}e_1$.
\enddemo

Indeed, the exactness in the middle term follows from the reciprocity in the same manner as in lemma 6.3. The exactness in the last term follows from
$|H_1(N)^{\Z_p}| = |H_1(N)_{\Z_p}| = |\text{Ker} \, \epsilon |$.
The rest of the argument of section 3 goes unchanged and we arrive to the following result.

\proclaim{7.2 Theorem (Splitting theorem for shrinking)}
The group
$H_1(N)$
splits as a direct orthogonal sum

$$H_1(N) = (W_1 + \ldots + W_{s-1}) \oplus W_s$$
where
$W_s \approx \Z/p^{k_s-1} \Z$
with the trivial action and
$W_i \approx \Lambda/\alpha_i$
for some ideals
$\alpha_i$.
\endproclaim

\demo{Remark 7.2}
The statement of the theorem holds true if
$p =2$
and
$H_1(M) = V \oplus \Z/2^k\Z, \, k \ge 2$,
with the same proof.
\enddemo

\head{8. Isotropic extension} \endhead

\demo{8.1}
Assume that
$H_1(M) = \Z/p^{k_1} \Z \oplus \ldots \oplus \Z/p^{k_s} \Z \oplus (\Z_p \oplus \Z_p)$
is the orthogonal decomposition and the form in the last summand is hyperbolic, given by the matrix

$$\pmatrix 0 & 1\\ 1 & 0 \endpmatrix .$$
We study the covering, corresponding to the map
$\epsilon : H_1(M) \rightarrow \Z_p$
defined by
$(\quad , e_{s+2})$
and, sending 
$e_i$
to zero,
$i \le s$
or
$i = s+2$
and
$e_{s+1}$
to a generator. Let
$v_i, \, i \le s$
be as above, let
$T = t (v_{s+1})$
and let $v$ be such that
$\pi_*v = e_{s+2}$
and
$N(v) = t \, e_{s+2}$,
such $v$ exists because
$e_{s+2} \in \text{Ker} \epsilon$.
For
$i \le s$
we keep on denoting
$p_i = te_i$.
\enddemo

\proclaim{Lemma (8.2)}
\roster
\item"1)" $ord \, p_i = p^{k_i}$
for
$i \le s$ \newline
\item"2)" $ord \, T = p$ \newline
\item"3)" $p_i$
and $T$ generate
$H_1(N)^{\Z_p}$
and
$\sum a_ip_i + aT =0$
implies
$p^{k_i}|a_i$
and
$p|a$.
\endroster
\endproclaim

\demo{Proof}
The proof of completely parallel to that of lemma 6.2. Let us prove for example, 3). Suppose
$\sum a_i p_i + aT =0$.
Then for
$j \le s$,
$0 = (v_j,\sum a_ip_i + aT)_N = \sum_i a_i(e_i,e_j) + a(e_j,e_{s+1}) = a_j(e_j,e_j)$,
so
$p^{k_i}|a_j$,
hence
$\sum a_ip_i =0$
and
$aT =0$,
and we deduce that also
$p|a$.
\hfill Q.E.D.
\enddemo

Let
$W_i,\, i\le s$
be the cyclic submodule, generated by
$v_i$.
Then
$(\sum W_i)^{\Z_p}$
contains a subgroup generated by
$p_i$,
and since
$|(\sum W_i)^{\Z_p}| = |(\sum W_i)_{\Z_p}| = \oplus_{i\le s} \Z/p^{k_i}\Z$,
$T \not \in \sum W_i$.
It follows, exactly as in lemma 6.4, that
$(\cdot , \cdot)|\sum W_i$
is nondegenerate. Let
$H_1(N) = (\sum W_i) \oplus (\sum W_i)^{\perp}$
be the orthogonal decomposition. Observe that
$T \in (\sum W_i)^{\perp}$
by reciprocity. Now, it follows easily, that one can modify $v$ to
$\tilde v$
such that
$\tilde v \in (\sum W_i)^{\perp}$.
Let $W$ be the cyclic submodule, generated by
$\tilde v$.
There should be an invariant element in $W$, and since by lemma 8.2. 3),
$p_i$
and $T$ generate
$H_1(N)^{\Z_p}$,
we conclude that
$T \in W$. \newline
We claim that
$N \tilde v = (1 + \zeta + \ldots + \zeta ^{p-1}) \tilde V =0$.
Indeed, for any
$w \in H_1(N), \, (N \tilde V, w)_N = (e_{s+2},\pi_*W)_N$
but
$\pi_*W \in \text{Ker} \epsilon$
and
$e_{s+2}$
is in the kernel of
$(\cdot , \cdot)_M|_{\text{Ker} \epsilon}$.
We sum up all the information in the following theorem

\proclaim{Theorem (8.3) (Splitting theorem for isotropic extensions)}
The group
$H_1(N)$
splits orthogonally as

$$H_1(N) = \sum_{i \le s} W_i \oplus W,$$
where $W$ is a
$\Lambda/\Lambda \cdot (1+\zeta + \ldots \zeta^{p-1}) = \Cal O$
-module.

Moreover,
$H^i(\Z_p,H_1(N)) \approx \Z_p$
for
$i \ge 1$
\endproclaim
\demo{Proof}
Only the last statement needs proof. It is enough to show that
$T \not \in Im \, N$
(recall that
$N = 1 + \zeta + \ldots + \zeta^{p-1})$.
But if
$T = Nw$,
then
$(T,v) = (Nw,v) = (w,Nv) =0$.
However
$(T,v) = (e_{s+1},e_{s+2})_M =1$,
a contradiction.
\enddemo

\head {9. Blowing up $H_1$: the non-abelian step} \endhead

In this chapter we prove that any ``generic'' 3-manifold $M$ has a finite covering $N$ with large
$H_1(N)_{(p)}$.
We first consider the case of homology sphere $M$ with the following ``genericity'' assumption:

\demo{(R)}
The Casson invariant
$|\lambda (M)| > \sharp$
(representation of
$\pi_1(M)$
onto
$SL_2(\F_5))$.
\enddemo

Our first result is:

\proclaim{Theorem (9.1)}
Let $M$ be a homology sphere satisfying 
$(R)$.
Then either $M$ is virtually Haken or
$\pi_1(M)$
admits a Zariski dense representation in
$SU(2)$.
\endproclaim

\demo{Proof}
Since
$\lambda (M) \ne 0$,
there are some nontrivial representations of
$\pi_1(M)$
in
$SU(2)$
[A]. Suppose none of them is Zariski dense. Since
$\pi_1(M)$
is perfect, we see that images of all representations should be finite. Moreover, among all finite subgroups of
$SU(2)$
only
$SL_2(\F_5)$
does not have abelian quotients [W]. Therefore, for any nontrivial representation
$\rho : \pi_1(M) \rightarrow SU(2), \, \rho (\pi_1(M)) \approx SL_2(\F_5)$.
Let
$W_1 \underset S \to \cup W_2$
be a Heegard splitting of $M$. Let
$R_1,R_2,R$
be representation varieties of
$\pi_1(W_1), \, \pi_1(W_2)$
and
$\pi_1(S)$
in
$SU(2)$.
Notice that set theoretically
$R_1 \cap R_2$
is the representation variety of
$\pi_1(M)$.
We claim:
\enddemo

\proclaim{Lemma (9.1)}
For any
$\rho \in R_1 \cap R_2, \, R_1$
and
$R_2$
intersect transversally at
$\rho$.
\endproclaim

\demo{Proof}
The Zariski tangent spaces of
$R_i$
at
$\rho$
are
$H^1(\pi_1(W_i),su(2))$
with adjoint action and since
$\pi_1(M) = \pi_1(W_1) \underset {\pi_1(S)} \to \ast \pi_1(W_2)$,
the twisted Myer Wietoris exact sequence of [JM] shows that
$T_{\rho}R_1 \cap T_{\rho}R_2 \approx H^1(\pi_1(M), su(2))$.
If the intersection is not transversal, then
$H^1(\pi_1(M),su(2)) \ne 0$.
Let $N$ be the covering of $M$ defined by the exact sequence

$$1 \rightarrow \pi_1(N) \rightarrow \pi_1(M) \overset {\rho} \to \longrightarrow SL_2(\F_5) \rightarrow 1,$$

then
$H^1(\pi_1(M),SU(2)) = (H^1(\pi_1(N),su(2))^{SL_2(\F_5)}$,
so that
$H^1(\pi_1(N),su(2)) \ne 0$.
But the action of
$\pi_1(N)$
in
$su(2)$
is trivial, and so that
$H^1(N,\R) \ne 0$,
a contradiction.
\enddemo

Returning to the proof of the theorem, we see that
$R_1$
and
$R_2$
intersect transversally in a finite number of points, each being a representation of
$\pi_1(M)$
to
$SL_2(\F_5)$,
so that
$|\lambda (M)| \le \sharp$
(representation of
$\pi_1(M)$
to
$SL_2(\F_5))$.
\hfill Q.E.D.

$   $\newline

An application of strong approximation in algebraic groups, as described in [Lu], gives immediately the following result

\proclaim{Theorem (9.2)}
Let $M$ be a homology sphere, satisfying
$(R)$.
Then for any prime $p$ there exists a finite normal covering
$N \rightarrow M$
with rank
$(H_1(N)_{(p)})$
arbitrarily large. In particular, there exists a covering with rank
$(H_1(N)_{(p)}) \ge 4$.
\endproclaim

We will now sketch an alternative approach, with slightly weaker statement, having in mind the application to manifolds without assumption $(R)$. First, we may assume that the complex Zariski dense representation
$\rho : \pi_1(M) \rightarrow SL_2(\C)$
is defined over
$\bar \Q$,
using the usual deformation argument [Re1]. So there is a subring
$\Cal O_S$
in a number field $F$, such that
$Im \, \rho \subset SL_2(\Cal O_S)$.
Now, the Bass-Serre theory of groups, acting on trees [Se4] implies immediately, with the argument of Culler-Shalen [CS], that either $M$ is vitually Haken, or
$\rho$
can be defined over
$\Cal O \subset \F$.
Let $p$ be a prime such that
$H_1(M)_{(p)} = 0$
and for some prime
$\frak p$
over $p$ in
$\Cal O, \, \Cal O/\frak p = \F_p$;
there are infinitely many such primes. We claim that the composite map
$\rho_p ; \pi_1(M) \rightarrow SL_2(\Cal O) \rightarrow SL_2(\F_p)$
is on. Indeed, since
$p \nmid |H_1(M)|$
and there are no proper perfect subgroups in
$SL_2(\F_p)$,
either
$\rho_p$
is on, or
$Im \, \rho_p =1$,
but the last option is impossible, since the kernel of reduction
$SL_2(\Cal O) \rightarrow SL_2(\F_p)$
is
$p$-residually finite. Now, we look at the covering
$N \rightarrow M$,
defined by the short exact sequence
$1 \rightarrow \pi_1(N) \rightarrow \pi_1(M) \overset {\rho_p} \to \longrightarrow SL_2(\F_p)$.
We assume
$p > 5$.
Our first claim is that
$H_1(N)_{(p)} \ne 0$.
Indeed, the s.s. of the covering
$N \rightarrow M$
reads like

$$\align &H^*(SL_2(\F_p),\Z) \\
&H^*(SL_2(\F_p),W) \\
& \qquad 0 \qquad \qquad \qquad \qquad \Rightarrow H^{i+j} (M,\Z) \\
& H^*(SL_2(\F_p),\Z)\endalign$$
where
$W = H^2(N) = \widehat{H_1} (N)$.
Now,
$H^*(SL_2(\F_5),\Z)_{(p)}$
is
$\F_p$
for
$* = k(p-1)$
and 0 otherwise, so
$W =0$
is impossible.

$   $\newline
Now we will look at the same s.s. with
$\F_p$
-coefficients:

$$\align
&H^*(SL_2(\F_p),\F_p) \\
&H^*(SL_2(\F_p),V^*) \\
&H^*(SL_2(\F_p),V) \qquad \qquad \qquad \Rightarrow H^{i+j} (M,\F_p) \\
&H^*(SL_2(\F_p),\F_p) \endalign$$
where
$V = H^1(N,\F_p)$.
Since
$p \nmid |H_1(M)|, \, V^{SL_2(\F_p)} =0$.
This eliminates the possibility of
$\dim V =1$.
Now, if 
$\dim V \le 3$,
then $V$ is either the natural two-dimensional module, or the adjoint module
$sl_2(\F_p)$.
In both cases, $V$ is cohomologically trivial for
$p > 3$
as the restriction to the $p$-Silow subgroup shows, which is impossible by the same argument as above. Hence
$\dim H_1(N,\F_p) \ge 4$
or $M$ is virtually Haken.

\head{10. Class towers of three-manifolds and the structure of the pro-$p$ completion of $\pi_1$}\endhead

\demo{10.1}
In this chapter we assume that $M$ is a three-manifold, not virtually Haken and such that rank
$H_1(M)_{(p)} \ge 4$
for some $p$. Our main result is as follows
\enddemo

\proclaim{Theorem (10.1)}
Let
$M_1 = M = \tilde M_1$
and let
$M_i$,
(resp.
$\tilde M_i\, ) i \ge 2$
be the maximal abelian $p$-covering (resp. maximal elementary abelian $p$-covering) of
$M_{i-1}$
(resp.
$\tilde M_{i-1}$).
Let
$r_i =$
rank
$(H_1(M_i)_{(p)})$
(resp.
$\tilde r_i =$
rank
$(H_1(\tilde M_i)_{(p)}$).
Then
\roster
\item"1)" $r_{i+1} \ge \frac {r^2_i-r_i} 2 $, 
and the same for
$\tilde r_i$, \newline
\item"2)" The group
$H_1(\tilde M_{i+1})_{(p)}$
has exponent
$\ge p^{\tilde r_i-1}$
\endroster
In particular, the tower
$\{M_i\}$
is infinite; that is, the pro-$p$ completion of
$\pi_1(M)$
is infinite.
\endproclaim

\demo{Proof}
Let
$H_1(M)_{(p)} = \Z/p^{k_1} \Z \oplus \ldots \oplus \Z/p^{k_r} \Z$
where
$r = r_1$.
Denote
$A = H_1(M)_{(p)}$
and consider the covering

$$M_1 = N \rightarrow M$$
with

$$1 \rightarrow \pi_1(N) \rightarrow \pi_1(M) \rightarrow A \rightarrow 1$$
exact.
\enddemo

We will write the s.s. of the covering
$N \rightarrow M$
in integer homology:

$$\CD
\Z  @. @. @. A @. \ldots \\
0   @. @.  @. 0 @. \ldots \\
(H_1(N))_A @. \qquad \qquad  @. H_1(A,H_1(N)) @. @. \ldots @. @. \Rightarrow H_{i+j} (M,\Z) \\
\Lambda^2_{\Z} A \\
\Z @. @. @.  A \endCD$$
we see immediately that
$d_2 : \Lambda^2_{\Z} A \rightarrow (H_1(N))_A$
should be injective (since
$H_2(M) =0)$.
So rank
$H_1(N)_{(p)} \ge$
rank
$(H_1(N)_{(p)})_A \ge \frac {r^2-r} 2$,
as stated. The proof for
$\tilde r_1$
is identical. \newline

$   $\newline
To prove (2) we will write the s.s. of the covering defined by the exact sequence
$1 \rightarrow \pi_1 (N) \rightarrow \pi_1(M) \rightarrow A \rightarrow 1$,
where now
$A = \undersetbrace{r} \to {\Z/p \oplus \ldots \oplus \Z/p \Z}$
in cohomology. We have

$$\CD
\Z @. \quad H^*(A,\Z) \\
H^*(A,H^2(N) \\
0 @. \qquad \qquad \qquad \Rightarrow H^{i+j} (M,\Z) \\
H^*(A,\Z) \endCD$$
By proposition 1.1, the exponent of
$H^*(A,\Z)$
is $p$. Since the
$E^{\infty}$
should look like

$$\CD
p^r\Z @. \qquad @. \quad 0 \\
\ast \\
0 \\
\Z @. \quad 0 @. \quad \ast
\endCD$$
we conclude that the exponent of
$H^2(A,H^2(N))$
is at least
$p^{r-1}$,
hence the exponent of
$H^2(N) \approx H_1(N)$
is at least
$p^{r-1}$,
which concludes the proof.

\demo{10.2}
The following result specifies the structure of pro-$p$ completion of
$\pi_1(M)$.
Informally speaking, a counterexample to the Thurston Conjecture would give rise to a ``$p$-adic three manifold'', with ``fundamental group''
$\frak G$.
\enddemo

\proclaim{Theorem (10.2)}
In conditions of 10.1, let
$\frak G$
be the pro-$p$ completion of
$\pi_1(M)$.
Then
$\frak G$
is a Poincar\'e duality pro-$p$ group in dimension 3. \endproclaim

\demo{Proof}
First we notice that the subsequent quotients of class towers

$$\align &\pi_1(M) = \pi_1(M) \supset \pi_1(M_2) \supset \ldots \\
\text{and} \quad &\frak G = \frak G_1 \supset \frak G_2 \supset \ldots \endalign$$
are identical. It follows that
$\frak G_k = [\frak G_{k-1}, \frak G_{k-1}]$
is a pro-$p$ completion of
$\pi_1(M_k)$.
Hence
$H_1(\frak G_k)$
is finite and
$H^1(\frak G_k,\Z) = 0$.
Let
$G_k = \frak G/\frak G_k$;
then
$G_k$
is a finite $p$-group. Now, for all
$k, \, H^2(G_k,\Z) \approx \widehat{H_1} (G_k) = \widehat{H_1} (M)_{(p)}$,
so
$H^2(\frak G,\Z) \approx \widehat{H_1} (M)_{(p)}$.
Consider the s.s. of the extension
$1 \rightarrow \frak G_k \rightarrow \frak G \rightarrow G_k \rightarrow 1$
in integral cohomology. We get

$$\matrix \format \c \quad & \c \quad & \c \quad & \c \quad &  \c \quad & \c \quad & \c \quad & \c \\
& \hat W^{G_k}_k & H^1(G_k,\hat W_k) & H^2(G_k,\hat W_k) & H^3(G_k,\hat W_k) \\
& 0 & 0 & 0 & 0 & 0 & \Rightarrow  H^{i+j}(\frak G,\Z) \\
& \Z & 0 & H^2(G_k,\Z) & H^3(G_k,\Z) & H^4(G_k,\Z)
\endmatrix$$
where
$\hat W_k = H^2(\frak G_k,\Z)$.
We wish to compare this s.s. of the covering 
$M_k \rightarrow M$

$$\matrix \format \c \quad & \c \quad & \c \quad & \c \quad &  \c \quad & \c \quad & \c \quad & \c \\
& H^2(M_k,\Z)^{G_k} & H^1(G_k,\hat W_k) & \ldots &\qquad \qquad & \Rightarrow &H^{i+j} (M,\Z) \\
& 0 & 0 & 0 & 0 & 0 \\
& \Z & 0 & H^2(G_k,\Z) & H^3(G_k,\Z) & H^4(G_k,\Z) 
\endmatrix$$
and
$H^2(M_k,\Z) \approx \widehat{H_1(M_k,\Z)} \approx \hat W_k$,
so that the first three rows of the two s.s are identical. Observe that there is a natural map from the first s.s. to the second so that the differential
$d_3$
from the third row to the first is the same for both.
\enddemo

Now, since the wedge map
$H^2(G_k,\Z) \rightarrow H^2(M,\Z)$
is an isomorphism dual to
$H_1(M,\Z) \approx G_k/[G_k,G_k]$,
we see that
$d_3 : H^2(M_k,\Z)^{G_k} \rightarrow H^3(G_k,\Z)$
is injective for the second s.s., and, in fact, an isomorphism, since the finite group
$coker \, d_3$
should inject to
$H^3(M,\Z) = \Z$,
so
$coker \, d_3 =0$.
That means in terms of the first s.s. that the wedge map
$H^3(G_k,\Z) \rightarrow H^3(\frak G,\Z)$
is zero, so
$H^3(\frak G,\Z) = \varinjlim H^3(G_k,\Z) =0$.
Next, the group
$coker \, d_3 : H^1(G_k,\hat W_k) \rightarrow H^4(G_k,\Z)$
from the second s.s. should be killed by
$d_4$
from the cyclic subgroup of
$H^3(M_k,\Z) = \Z$,
so it is cyclic. In terms of the first s.s. it means that the image of
$H^4(G_k,\Z)$
in
$H^4(\frak G,\Z)$
is cyclic, so
$H^4(\frak G,\Z)$
is either
$\Z/p^N\Z$
or
$(\Q/\Z)_{(p)}$.
We claim that the first option is impossible. Indeed, if this is the case then
$|\text{coker} \, d_3 : H^1(G_k,\hat W_k) \rightarrow H^4(G,\Z)| \le p^N$.
(since
$H^3(\frak G_k,\Z) =0$,
the forth row of the first s.s. is zero). On the other hand, since the image of
$H^3(M,\Z)$
in
$H^3(M_k,\Z)$
has index
$|G_k|$,
the product of orders of groups
$\text{coker} \, d_3 : H^1(G_k,\hat W_k) \rightarrow H^4(G,\Z)$
and
$Im \, d_2 : H^3(M_k,\Z) \rightarrow H^2(G_k,\hat W_k)$
is
$|G_k|$,
so the former group has order
$\ge \frac {|G_k|} {p^N}$
since these groups are recipients of differentials from
$H^3(M_k,\Z)$.

$   $\newline

So
$|G_k| \le p^{2N}$
which is impossible by Theorem 10.1. So
$H^4(\frak G,\Z) = (Q/\Z)_{(p)}$. \newline

$   $\newline
We claim that for
$i \ge 5, \, H^i(\frak G,\Z) =0$.
Indeed, the first s.s. above now looks like

$$\matrix \format \c \quad & \c \quad & \c \quad & \c \quad &  \c \quad & \c \quad & \c \quad & \c \\
\vdots \\
H^*(G_k,(\Q/\Z)_{(p)}) \\
0 \\
H^*(G_k,\hat W_k) & \Rightarrow & H^{i+j}(\frak G,\Z) \\
0 \\
H^*(G_k,\Z)
\endmatrix$$
and
$d_3 : H^*(G_k,\Q/\Z)_{(p)}) \rightarrow H^{*+3} (G_k, \hat W_k)$
is the same that
$d_2 : H^* (G_k,\Z) \rightarrow H^{*+2} (G_k, \hat W_k)$
from the second s.s. after identification
$H^*(G_k(\Q/\Z_{(p)}) \approx H^{*+1} (G_k,\Z)$.
This implies immediately that
$H^i(G_k,\Z), \, i \ge 5$
is killed in the first s.s. because it is killed in the second. Summing up, we have the following table:

$$\matrix \format \c \quad & \c \quad & \c \quad & \c \quad &  \c \quad & \c \quad & \c \quad & \c \\
H^0(\frak G,\Z) &H^1(\frak G,\Z) & H^2(\frak G,\Z) & H^3(\frak G,\Z) & H^4(\frak G,\Z) & 0 & \ldots \\
\Z & 0 & \widehat{H_1} (M)_{(p)} & 0 & (\Q/\Z)_{(p)} & 0 & \ldots \endmatrix$$
hence in
$\F_p$
-coefficients we have

$$\matrix \format \c \quad & \c \quad & \c \quad & \c \quad &  \c \quad & \c \quad & \c \quad & \c \\
H^0(\frak G,\F_p) & H^1(\frak G,\F_p) & H^2(\frak G,\F_p) & H^3(\frak G,\F_p) & 0 & \ldots \\
\F_p & V & V^* & \F_p & 0 & \ldots
\endmatrix$$
where
$V = Hom (H_1(M),\Z/p\Z)$
and it is easy to prove that the duality
$V \times V^* \rightarrow \F_p$
is induced by multiplication. \newline
\hfill Q.E.D.

\head{11. Nazarova-Royter theory and the structure of anisotropic and isotropic extensions, II} \endhead

\demo{11.1}
We should recall the Nazarova-Royter classification of finite modules over the cyclic group
$C_p$,
c.f. [NR], [Le].
\enddemo

\demo{Definition}
A left elementary block is
$\Z/p^N\Z$
with trivial action of
$C_p$.
\enddemo

\demo{Definition}
A right elementary block is a module of the form
$\Cal O /(1-\zeta)^N \Cal O$,
where
$\Cal O = \Z[C_p]/1+\zeta + \ldots + \zeta^{p-1}$
the ring of cyclotomic integers. \newline
The prime
$1-\zeta \in \Cal O$
is called
$\pi$.
\enddemo

\demo{Definition}
A diagram

\midspace{3cm}
$   $\newline
represents a fibered sum of a left module $L$ and a right module $R$ over
$L/pL \approx R/\pi R = \F_p :$

$$0 \rightarrow W \rightarrow L \oplus R \rightarrow \F_p \rightarrow 0$$
\enddemo

\demo{Definition}
A diagram

\midspace{3cm} 
$   $\newline
represents a quotient of
$L \oplus R$
which identifies
$\text{Ker} \, (\times p)$
in $L$ with
$\text{Ker} ( \times \pi)$
in $R$.
\enddemo

\demo{Definition}
An open module is represented by a diagram

\midspace{3cm}
$   $\newline
$(p \, x_i =0)$.
\enddemo

\demo{Definition}
A closed module is a quotient of an open module $W$ by a relation
$\sum a_ix_i =0,\, a_i \in \F_p$,
where
$\sum a_it^i$
is a power of irreducible pdinomial
$\ne ax$
over
$\F_p$.
\enddemo

\proclaim{Theorem (11.2) (Nazarova-Royter)}
Any f.g. indecomposable module over
$C_p$
which is a $p$-torsion is either open or closed module.
\endproclaim

\proclaim{Theorem (11.3) ([AL])}
If $W$ is an open module. Then
$H^i(C_p,W) =\Z/p\Z$
for
$i \ge 1$.
If $W$ is a closed module, then
$H^i(C_p,W) =0$
for
$i \ge 1$.
\endproclaim

A combination with Proposition 2.1 yields immediatley the following result of key importance:

\proclaim{Theorem (11.4) (The structure of anisotropic and isotropic extensions, II)}
Let
$1 \rightarrow \pi_1(N) \rightarrow \pi_1(M) \rightarrow C_p \rightarrow 1$
is a normal covering of rational homology spheres, represented by an element $z$ of order $p$ in
$H_1(M)$,
that is, the map
$H_1(M) \rightarrow C_p$
is given by
$(\cdot, z)_M$.
Then the structure of
$H_1(N)_{(p)}$
as a 
$C_p$
-module is as follows:
\roster
\item"(a)" if
$(z,z) =0$,
then
$H_1(N)_{(p)}$
is a direct sum of open modules \newline
\item"(b)" if
$(z,z) \ne 0$,
then
$H_1(N)_{(p)}$
is a direct sum of exactly one closed module and several open modules.
\endroster
\endproclaim

\head{12. Non-existence of split anisotropic constituents} \endhead

\demo{12.1} 
Let $G$ be a finite group acting on a manifold $N$ freely. Suppose
$V \subset H_1(N)$
be a $G$-invariant subspace, and put
$W =\frac {H_1(N)} V$.
Let
$Q \rightarrow N$
be a covering, defined by a map
$\pi_1(N) \rightarrow H_1(N) \rightarrow W$.
Then
$Q \rightarrow M$
is a normal covering with a Galois group 
$\Cal P$, 
which is an extension

$$0 \rightarrow W \rightarrow \Cal P \rightarrow G \rightarrow 1.$$
\enddemo

\demo{12.2}
Now suppose
$1 \rightarrow \pi_1(N) \rightarrow \pi_1(M) \rightarrow C_p \rightarrow 1$
is a covering of rational homology three-spheres, and $N$ is virtually non-Haken (in particular, has virtual first Betti number zero).
\enddemo

The main result of this section is as follows.

\proclaim{Theorem (12.2)}
(a) If
$p \ne 2$,
then
$H_1(N)_{(p)}$
does not have a split anisotropic constituent, that is, an invariant cyclic subgroup generated by an element $z$ with
$ord \, (z,z)_N = ord \, z$. \newline
(b) Suppose
$p=2$
and
$H_1(N)_{(2)}$
has such a subgroup $W$. Consider a
$C_2$
-invariant orthogonal splitting
$H_1(N)_{(2)} = W \oplus V$.
Let
$Q \rightarrow N$
be as in 12.1. Then the group
$\Cal P$ is a generalized quaternionic (binary dihedral) group, and
$H^i(\Cal P, H_1(Q)) = 0$
for
$i \ge 1$.
In particular, the action on $W$ is the multiplication by
$(-1)$.
\endproclaim

\demo{Proof}
We will prove only the case (b), that is,
$p=2$,
since the case of odd $p$ is completely similar (and easier).
\enddemo

\proclaim{Lemma (12.2)}
In condition of Theorem 12.2 (b) we have:

$$\align &H^i(W,Q) =0 \quad \text{for} \quad i \ge 1 \\
&H^0(W,Q) \approx V \quad \text{as} \quad C_2-\text{module}. \endalign$$
\endproclaim

\demo{Proof}
This follows by induction from the Shrinking Theorem 7.2 and the argument of 3.1.
\enddemo

\proclaim{Lemma (12.3)}
$H^i(\Cal P, H_1(Q))$
is either 0 or
$\Z/2\Z$
for all
$i \ge 1$.
\endproclaim

\demo{Proof}
The Lyndon-Serre-Hochschild s.s. for the
$\Cal P$
-module
$H_1(Q)$
degenerates to
$H^i(C_2,V)$
by the previous lemma. Since $V$ is a split component of
$H_1(N)$,
the statement follows from Proposition 2.1. \newline
Now, we have exactly three possibilities for
$\Cal P$,
except that of quaternionic group:
\roster
\item"(i)"
$\Cal P = C_2 \times C_{2^n}$ \newline
\item"(ii)"
$\Cal P$ 
is dihedral 
$D_{2^n}$ \newline
\item"(iii)"
$\Cal P = C_{2^n} \rtimes C_2$
with the
$C_2$
-action on
$C_{2^n}$
given by multiplication by
$\pm (2^{n-1} -1)$
\endroster
\enddemo

In cases (i) and (ii) the cohomology
$H^i(C_2,W) \approx \Z/2\Z$
and
the proof of the Lemma 12.3 shows that
$H^i(\Cal P,H_1(Q)) = H^i(C_2,V) =0$
for
$i \ge 1$.
Hence the s.s. of the group extension shows that
$H^i(\Cal P,\Z)$
should be 4-periodic, which is not the case. \newline
In case (iii) a theorem of Wall [W], see also [Th], shows that
$H^{even}(\Cal P,\Z)$
is freely multiplicatively generated by two elements
$\xi,\eta$
of degree 2 subject to relations
$2\xi =0$
and
$2^n \eta =0$. \newline
Since the ranks of
$H^{2k}(\Cal P,\Z)$
are not bounded, there should be a nontrivial
$d_3 : (coker \, d_2)^{odd} \rightarrow H^{even}(\Cal P,\Z)$
in the s.s. of the covering
$Q \rightarrow M$
so that its image contains some
$0 \ne v \in H^{even}(\Cal P,\Z)$.
But then the Proposition 1.4 implies that the ranks of
$coker \, d_2$,
hence of
$H^*(\Cal P,H^2(Q))$
are unbounded because
$d_3(coker \, d_2)$
contains all
$\zeta^k \cdot \eta^l \cdot v$,
which is a contradiction to Lemma 12.3. 

We are now ready to determine the pro-$p$ completion of
$\pi_1(M)$
when
$H_1(M)_{(p)}$
is ``small''.

\proclaim{Theorem (12.4)}
Let
$H_1(M)_{(p)}$
be cyclic. Then either $M$ is virtually non-Haken, or $M$ is a quotient of $p$-homology sphere by a cyclic group. The pro-$p$ completion
$\frak G$
of
$\pi_1(M)$
is then finite cyclic.
\endproclaim

\demo{Proof}
Follows immediately from the Shrinking Theorem 7.2.
\enddemo

\proclaim{Theorem (12.5)}
Let
$H_1(M)_{(2)} = \Z/2\Z \oplus \Z/2\Z$.
Then: \newline
(i) if the linking form in
$H_1(M)_{(2)}$
is
$\langle 1 \rangle \oplus \langle 1 \rangle$
then is virtually Haken or $M$ is a quotient of a 2-homology sphere by a generalized quaternionic group
$Q_{2^n}, \, n \ge 4$,
and
$\frak G \approx Q_{2^n}$. \newline
(ii) if the linking form in
$H_1(M)_{(2)}$
is hyperbolic
$\pmatrix 0 & 1 \\1 & 0 \endpmatrix$
then either $M$ is virtually Haken, or $M$ is a quotient of a 2-homology sphere by the quaternionic group
$Q_8$
and
$\frak G \approx Q_8$.
\endproclaim

\demo{Proof}
Let
$z \in H_1(M)$
be any isotropic element and consider a
$C_2$
-covering
$N \rightarrow M$
defined by
$(\cdot , z) : H_1(M) \rightarrow \Z/2\Z$.
Then
$Q = H_1(N)_{(2)}$
is a 
$C_2$
-module and we know that it is a sum of exactly one open module and some closed modules. Since
$W^{C_2} \approx \Z/2 \Z$,
in fact $W$ is just one open module. Now, since the linking form in $W$ is
$C_2$
-invariant we have
$W \approx \hat W$,
which gives, along with the main result of [NR] and [Le], that $W$ is either of a form

\midspace{3cm}
$   $\newline
or

\midspace{3cm}
$   $\newline
Now, the condition
$W^{C_2} \approx \Z/2\Z$
implies that $W$ is actually cyclic with either trivial action ($L$-module) or 
$(-1)$
-action ($R$-module). By the Theorem 12.4 we have
$N = Q/W$
where $Q$ is a 2-homology sphere, and by the Theorem 12.2,
$\frak G$
is a (generalized) quaternionic group
$Q_{2^n}, \, n \ge 3$. \newline
Now, if the linking form is hyperbolic, then \underbar{any} element is isotropic, hence whatever homomorphism
$\frak G \twoheadrightarrow \Z/2\Z$
we consider, the kernel would be cyclic. This rules out the possibility
$n \ge 4$,
so
$\frak G \approx Q_8$.
Conversely, if
$n = 3$,
there is a map of odd degree from $M$ to
$S^3/Q_8$,
which induces an isomorphism in 2-torsion of
$H_1$.
Since
$S^3/Q_8$
has a selfhomeomorphism
$\zeta$
of order 3 (the quotient of the normalizer
$N(Q_8)$
by
$Q_8$
in
$SL_2(\F_5)$
has order 3), which does not have fixed nonzero elements in
$H_1(S^3/Q_8)_{(2)}$,
the linking form
of
$S^3/Q_8$
is hyperbolic.
\enddemo
\hfill Q.E.D.

\head{13. A spectral sequence in group cohomology} \endhead

In this section, we introduce a new powerful spectral sequence, converging to
$\F_p$
-cohomology of a normal subgroup of index $p$ in a given group $G$. There are ``commutative analogs'' of this s.s. in algebraic geometry, which appear to be well-known [Se1] \footnote{I am grateful to David Eisenbud for this reference}. \newline
Let
$1 \rightarrow K \rightarrow G \overset s \to \longrightarrow C_p \rightarrow 1$
be an exact sequence of groups.

\proclaim{Theorem (13.1)}
There is a spectral sequence with
$E_1$-term

$$E^{i,j}_1 = H^{i+j}(G,\F_p) \Rightarrow H^{i+j}(K,\F_p), \, i+j \ge 0, \, 0 \le i \le p-1$$
The differential
$d_1 : H^i(G,\F_p) \rightarrow H^{i+1} (G,\F_p)$
is a multiplication by the class $s$ viewed as an element of
$H^1(G,C_p)$.
\endproclaim

\proclaim{Theorem (13.2)}
Let
$p=2$.
There is an exact sequence

$$\ldots H^i(G,\F_2) \overset r \to \longrightarrow H^i(K,\F_2) \overset t \to \longrightarrow H^i(G,\F_2) \overset {\times s} \to \longrightarrow H^{i+1} (G,\F_2) \rightarrow \ldots$$
where $r$ is the restriction and $t$ is the transfer map.
\endproclaim

\demo{Proof}
We begin with a lemma.
\enddemo

\proclaim{Lemma (13.3)}
Consider the regular module
$\Z/p^N\Z[C_p]$
over
$C_p$.
There is a canonical filtration

$$0 = V_0 \subset V_1 \subset \ldots V_{N+1} = \Z/p^N \Z[C_p]$$
of
$C_p$
-modules with all successive factors
$\F_p$.
\endproclaim

\demo{Proof}
Let
$\zeta$
be a fixed generator of
$C_p$.
Consider the exact sequence

$$0 \rightarrow \Z/p^N\Z \overset {1+\zeta + \ldots \zeta^{p-1}} \to \longrightarrow\Z/p^N\Z [C_p] \rightarrow \Z/p^N\Z [C_p]/(1 + \zeta \ldots \zeta^{p-1}) \rightarrow 0$$
Now, the latter module may be represented as
$(\Z[C_p]/(1+\zeta + \ldots + \zeta^{p-1}))/(p^N \Z[C_p]/(1+\zeta + \ldots +\zeta^{p-1}))$.
Let
$\Cal O = \Z[C_p]/(1+\zeta + \ldots + \zeta^{p-1})$
be the ring of cyclotomic integers. Since
$p = unit \times (1-\zeta)^{p-1}$
we see that the module above is just
$\Cal O/(1-\zeta)^{N(p-1)}$
and has canonical filtration by
$(1-\zeta)^l \Cal O/(1-\zeta)^{N(p-1)}$.
\enddemo
\hfill Q.E.D.

\demo{Proof of Theorem 13.1}
Since
$H^i(K,\F_p) = H^i(G,\F_p[G/K])$
by Shapiro's lemma, we may deal with the latter group. Now, we may rewrite it as
$Ext^i_{\Z G} (\Z,\F_p[G/K])$.
The filtration of
$G/K$
modules
$0 = V_0 \subset V_1 \subset \ldots \subset V_p = \F_p[G/K]$
becomes a filtration of $G$-modules, so the standard s.s. of filtrated modules [CE] gives a s.s. with
$E_1$
-term
$Ext^{i+j}_{\Z_ G}(\Z,V_i/V_{i+1}) \Rightarrow Ext^{i+j}_{\Z G} (\Z,\F_p[G/K])$,
as stated.
\enddemo

\demo{Proof of Theorem 13.2}
This follows from the previous proof, or may be seen from the short exact sequence of $G$-modules

$$0 \rightarrow \F_2 \overset {\times (1+\zeta)} \to \longrightarrow \F_2[G/K] \overset {aug} \to \longrightarrow \F_2 \rightarrow 0$$
\enddemo

\demo{Corollary 13.3}
Let $G$ be a pro-$p$ group and $K$ a subgroup of index $p$. Then there is a s.s. as in Theorem 12.1, converging to
$H^*(K,\F_p)$.
\enddemo

\demo{Proof}
Inductive limit by the directed set of finite index subgroups.
\enddemo

\head{14. Strenghtened Adem inequalitites and other applications} \endhead

Recently Adem [Ad] published several strong results, showing a possible range for finite $p$-group cohomology. All these results may be derived from our spectral sequence, as we will see soon. For
$p = 2$
we show some simple consequences for the structure of the cohomology ring; these are related to a well-known theorem of Serre, see [Se2], [Se3].

\proclaim{Theorem (14.1)}
Let $G$ be a finite $p$-group and let
$K < G$
be a subgroup of index $p$. Then
\roster
\item"(i)" $b_i(K,\F_p) \le p \, b_i(G,\F_p), \, i \ge 1$ \newline
\item"(ii)" $b_i(G,\F_p) > 0$
for all $i$ \newline
\item"(iii)" if
$p=2$
and
$G/[G,G]$
is elementary abelian, then
$b_1(K,\F_p) \le p(b_1(G,\F_p)-1)$ \newline
\item"(iv)" let
$r_i(G) = \log_p |H^i(G,\Z)|$.
Then for $i$ even,

$$p \cdot r_i(G) \ge r_i(K) + 1 ,$$
in particular
$r_i(G) > 0$
for all even $i$; for $i$ odd,

$$p \cdot r_i(G) \ge r_i(K) -1 .$$
\endroster
\endproclaim

\demo{Proof}
(i) is immediate from Theorem 13.1. (ii) follows from (i) by induction, since it is yielded for a cyclic group. For proving (iii), we first notice that
$d_1 : H^0 \rightarrow H^2$
is nonzero, which kills
$(p-1)$
independent elements in $\qquad$. For
$p =2, \, s\cdot s = s^2 =p(s) \ne 0$,
so
$d_1 : H^1 \rightarrow H^2$
also kills $s$ from the second
$H^1(G)$
in 13.2, hence
$b_1(K) \le 2(b_1(G)-1)$.
To prove (iv), we will first show the idea, proving it for
$p=2$. Consider the truncated exact sequence

$$C : 0 \rightarrow \F_p \rightarrow H^1(G) \rightarrow H^1(F) \rightarrow H^1(G) \rightarrow \ldots \rightarrow H^i(G) \rightarrow H^i(F) \overset {\psi} \to \longrightarrow Im \, \psi \subset H^i(G) \rightarrow 0$$
We have:

$$0 = \chi (C) \ge 1 -(b_1(G) -b_2(G) + \ldots +(-1)^{i+1}b_i(G)) + b_1(F) + \ldots + (-1)^{i+1} b_i(F)$$
for $i$ even and the opposite sign for $i$ odd.

Now, because of the short coefficient sequence
$0 \rightarrow \Z \overset{\times p} \to \longrightarrow \Z \rightarrow \F_p \rightarrow 0$
we have
$b_i = r_i + r_{i+1}$
and
$r_1 = 0$,
which gives (iv) for
$p =2$.
\enddemo

Recall the result of Serre [Se2], [Se3]: if $G$ is not elementary abelian then there are nonzero elements
$x_1 \ldots x_n \in H^1(G)$,
such that
$\beta x_1 \ldots \beta x_n =0$.
For
$p=2$ we have the following

\proclaim{Theorem (14.2)}
Suppose $G$ is not elementary abelian and let $r$ be the elementary abelian rank of $G$. Then
\roster
\item"(i)" The commutative 
$\F_2$
-algebra
$H^*(G,\F_2)$
has divisors of zero \newline
\item"(ii)" Moreover, for
$s \in H^1(G)$
let
$k_i = \dim \ker (H^i(G) \overset {\times s} \to \longrightarrow H^{i+1} (G))$.
Then there exists
$s \in H^1(G)$
such that

$$k_i + k_{i+1} \ge const. \, \cdot r^i, \quad \text{for} \quad i\ge 1.$$
\endroster
\endproclaim

\demo{Proof}
Let
$A \underset {\ne} \to \subset G$
be a maximal elementary abelian group. Choose
$s : G \rightarrow \Z/2\Z$
such that
$A \subset Ker \, s$
and put
$K = Ker \, s$.
Then
$r(k) = r(G) = r$.
From Theorem 13.2 we have

$$b_i(K) = k_i + b_i(G) -(b_{i-1}(G) -k_{i-1}) = k_i+k_{i-1} + (b_i(G) -b_{i-1} (G))$$
Now, from the main result of Quillen [Q] it follows that
$b_i(K) \sim const. \, \cdot r^i$
and
$b_i(G) -b_{i-1}(G) \sim const. \, \cdot r^{i-1}$,
which proves (ii), hence (i).
\enddemo

\head{15. Multiplication in cohomology: $\Z/2 \oplus \Z/2\Z \oplus \Z/2\Z$ case} \endhead

In this chapter we start a deeper analysis of homology of coverings for manifolds with small
$H_1$.
The main result of this chapter is a next theorem, identifying the multiplication in cohomology in case
$H_1(M)_{(2)} = (\Z/2\Z)^3$.
Observe that in
$(\Z/2\Z)^3$
any nondegenerate linking form is diagonalizable.

\proclaim{Theorem 15.4}
Let $M$ be a non-virtually Haken manifold with
$H_1(M)_{(2)} = \Z/2\Z \oplus \Z/2 \Z \oplus \Z/2\Z$.
Then there exists a basis
$x,y,z$
in
$H^1(M,\F_2)$
such that either
$xy = xz = yz =0$
or
$x^2=yz, \, y^2 =xz, \, z^2 =xy$.
\endproclaim

\demo{Proof}
Let
$x,y,z$
be an orthonormal base for
$(\cdot , \cdot)_M$.
We always identify
$v \in H_1(M)_{(2)}$
and
$(\cdot , v) \in H^1(M,\F_2)$.
Let
$w = x+y$,
then $w$ is isotropic. Let
$N \rightarrow M$
be a 
$C_2$
-covering, defined by $w$. Let
$W = H_1(N)_{(2)}$.
Then $W$ is a sum of one open module and several closed modules by the Theorem 11.4, and moreover
$W^{C_2} \approx \Z/2 \Z \oplus \Z/2 \Z$.
That means that either we have one open module, or a direct sum of an open module and a closed module.
\enddemo

\demo{Case A}
$W$ is an open module, hence of the type

\midspace{2cm}
$   $\newline
or

\midspace{2cm}
$   $\newline
by the argument of 12.5. The first case is definitely impossible, since then
$W_{C_2}$
would have direct summand
$\Z/2^n\Z , \, n \ge 2$.
(Recall that
$W_{C_2} \approx W^{C_2}$
since
$W \approx \hat W$).
In the second case we see immediately that the only possibility is

\midspace{2cm}
$   $\newline
and as an abelian group
$W \approx \Z/2^n \oplus \Z/2^n \Z$.
\enddemo

\demo{Lemma (15.5)}
$z^2 \ne w \cdot v$
in
$H^2(M,\F_2)$.
\enddemo

\demo{Proof}
Suppose
$z^2 = (x+y)(ax + by + cz)$.
Multiplying by $x$ and accounting the identity
$p^2 \cdot q = p \cdot q^2 =(p,q)_M$
we have
$0 = xz^2 = a + c \, xyz$.
Multiplying by $y$ we have similarly
$0 = b + c\, xyz$,
hence
$a = b$.
Now, multiplying by $z$ we have
$1 = (z,z) = z^3 + (a+b)xyz = 2axyz =0$,
a contradiction.
\enddemo

\demo{Corollary to lemma}
$z^2$
restricts nontrivially to
$H^2(N,\F_2)$.
\enddemo

\demo{Proof}
This follows from the lemma and the Theorem 13.2. 

Now, let
$\bar z$
be a restriction of $z$ to
$H^1(N,\F_2)$.
We have
$\beta \bar z = \bar z^2 \ne 0$,
which means that $W$ has a 
$\Z/2\Z$
-direct summand, a contradiction. So the case A is impossible.
\enddemo

\demo{Case B}
$W = W_1 \oplus W_2$,
where
$W_1$
is open and
$W_2$ is closed. Again, since
$W \approx \hat W, \, W_1$
should be as shown in the beginning of the analysis of case A. Since
$W^{C_2}_1 \approx \Z/2\Z$,
it follows immediately that
$W_1$
is cyclic (either left or right). Similarly, since
$W^{C_2}_2 \approx \Z/2\Z , \, W_2$
should look like

\midspace{1.5cm} 
$   $\newline
and necessarily 
$n=2$.
If
$m > 2$,
then
$W_2$
is cyclic with the acion being multiplication by
$\pm (2^{m-1}-1)$.
Since by the argument in the case A, $W$ contains a
$\Z/2\Z$
-summand, we see that
$W_1 = \Z/2\Z$
and
$W_2  =\Z/2^m\Z$,
but then 
$W_2$
is split
anisotropic with action different from
$(-1)$,
in contradiction to Theorem 12.2 (b). So
$m = n= 2$
and
$W_2 \approx \Z/2 \Z [C_2]$.
In short, $W$ looks like

\midspace{2cm}
$   $\newline
in particular
$b_1(N) =3$.
Again by the Theorem 13.2 we have that the kernel
$H^1(M,\F_2) \overset {\times w} \to \longrightarrow H^2(M,\F_2)$
is one-dimensional, say
$\alpha x + \beta y + \gamma z$.
We have
$0 = (x+y)(\alpha x + \beta y + \gamma z)$.
Arguing as above, we find
$0 =\alpha + \gamma \, xyz; \, 0 = \beta + \alpha \, xyz$,
so
$\alpha = \beta$.
If
$\gamma =0$
that gives
$\alpha = \beta =0$,
which is impossible, so
$\gamma =1$
and
$\alpha = \beta = xyz$.
\enddemo

\demo{Case $B_1$}
$xyz = \alpha = \beta =1$,
so that
$x^2 + y^2 = xz + yz$.
Then replacing $w$ by
$x + z$
or
$y+z$
and accounting that
$xyz =1$
we find also
$x^2+z^2 = xy + yz, \, y^2 + z^2 = xy +xz$.

Now, we look at the anisotropic covering
$N \rightarrow M$,
defined by $x$. Again we put
$W = H_1(N)_{(2)}$;
it is a sum of closed modules. What is a kernel of
$H^1 \overset {\times x} \to \longrightarrow H^2$ ?
Let
$x(ax + by + cz) =0$,
then we find as above
$a =0$.
Now,
$x(y+z) =x^2 + y^2$.
So
the only possibilities are
$xy =0$
or
$xz =0$,
Let us say, the first. Replacing $x$ by $z$ and repeating the procedure, we find that either the kernel
$H^1 \overset{\times z} \to \longrightarrow H^2$
is zero, or
$yz =0$,
or
$xz =0$.
But that contradicts the identities above. So we can assume that
$\ker H^1 \overset{\times x} \to \longrightarrow H^2$
is zero and
$\ker H^1 \overset {\times z} \to \longrightarrow H^2$
is zero. Now, we have
$xz = \alpha x^2 + \beta y^2 + \gamma z^2$
for some
$\alpha, \beta, \gamma$.
Equalities
$xz = \alpha x^2$
and
$xz = \gamma z^2$
would give
$x(z-\alpha x)=0$,
resp.
$z(x-\gamma z)=0$,
which is impossible by the above. Moreover,
$xz = x^2 +z^2$
would give
$xy = x^2 + y^2$
and
$y^2 = y^2 + z^2$.
Otherwise,
$\beta =1$,
so
$xz = \alpha x^2 + y^2 + \gamma z^2$.
If
$\alpha = \gamma =1$
we have
$xy = y^2+ z^2-(x^2+ y^2 + z^2) = x^2$,
a contradiction. Suppose
$\alpha =0, \, \gamma =1$,
so that
$xz = y^2+z^2$,
then
$xy =0$,
a contradiction. Likewise
$\alpha =1, \, \gamma =0$
is impossible. Hence
$\alpha = \gamma =0$
and
$xz = y^2, \, xy = z^2, \, yz = x^2$.

We claim that the case
$xz = x^2 + z^2$
above is impossible. Indeed, in that case there would be a map
$\pi_1(M) \rightarrow Q_8$,
such that the induced map in
$H^1(\cdot \, , \F_2)$
would send
$H^1(Q_8,\F_2)$
on the subspace spanned by
$x,y$
(see lemma 15.5 below). But this contradicts
$x^3 =1$.

So we conclude that
$x^2= yz, \, y^2 = x2, \, z^2 = xy$.
Moreover, since
$(x+y)^2 + (x+z)^2 + (x+y)(x+z) =0$,
we have a map
$M \rightarrow S^2/Q_8$,
inducing a homomorphism
$H^1(Q_8,\F_2) \rightarrow H^1(M,\F_2)$
with
$x+y$
and
$x+z$
in the image. This map has odd degree, since
$(x+y)^2 \cdot (x+z) =1$.
In particular,
$W_1$
above is
$\Z/4\Z$. 

Now, the exact sequence of the Theorem 13.2 gives
$b_1(N) =2$,
so $W$ as an abelian group has two generators. We have two cases again
\enddemo

\demo {Case $B_{11}$}
$W$ is a sum of two closed modules,
$W = W_1 \oplus W_2$
and both
$W_i$
are cyclic groups necessarily of degree
$\ge 8$. If their order is different, then the one with a bigger order is anisotropic, which is impossible by the Theorem 12.2 (b). So
$|W_1| = |W_2| \ge 8$.
Moreover, if
$W_1$
is not anisotropic, then
$W_1 \approx \widehat {W_2}$
and the action in
$W_1$
and
$W_2$
is he same, that is multiplication by either
$2^{n+1} +1$,
or
$2^{n-1}-1$
(where
$2^n = |W_i|)$.
So the action in $W$ is multiplication by
$2^{n-1} \pm 1$.
Now, we can choose
$W_i$
such that
$W_1$
projects on $y$ and
$W_2$
projects on $z$. By the argument of \newline
$\qquad \qquad$, we have a map
$\pi_1(M) \twoheadrightarrow \Phi$,
where
$\Phi$
is a semidirect product
$W_1 \rtimes \Z/2\Z$,
which induces the map
$H_1(M) \rightarrow \Z/2 \Z \oplus \Z/2 \Z \oplus \Z/2 \Z$,
sending $z$ to zero.
\enddemo

\demo{Lemma}
$b_2(\Phi) = \dim H^2(\Phi ,\F_2) =2$.
\enddemo

\demo {Proof}
We first notice that
$H^3(\Phi,\Z )=0$,
because
$H^i(\Z/2\Z, W_1)=0$
and hence
$H^i(\Z/2\Z , H^2(W_1)) =0$
for
$i \ge 1$.
Now from the exact coefficient system
$0 \rightarrow \Z \overset{\times 2} \to \longrightarrow \Z \rightarrow \Z/2 \Z \rightarrow 0$
we get
$b_2(\Phi ) =$
rank
$H^2(\Phi,\Z )=$
rank
$H_1(\Phi,\Z) = 2$.

By the lemma, there should be a relation between $x$ and $y$ in
$H^2(M,\F_2)$.
However,
$xy = z^2$
and
$x^2,y^2,z^2$
are independent, a contradiction. So the case
$B_{11}$
is impossible.
\enddemo

\demo{Case $B_{12}$}
$W$ is one closed module. Since
$W^{C_2} \approx \Z/2\Z \oplus \Z/2 \Z$,
$W$ should of a form

\midspace{2cm}
$   $\newline
but this has more than two generators as an abelian group, a contradiction.
\enddemo

\demo{Case $B_2$}
$xyz =0$
and
$(x+y)z =0$.
Then similarly
$(xz)y =0$
so
$xz = yz = xy$.
But
$(xz)y = (yz)x = (xy)z =0$,
so the element
$xz = yz = xy$
should be zero by Poincar\'e duality. This proves the theorem.
\enddemo

\demo{Corollary}
There exists homomorphisms
$\pi_1(M) \twoheadrightarrow Q_{2^n}, \, n \ge 4$,
such that the induced map
$H_1(M) \rightarrow \Z/2\Z \oplus \Z/2\Z$
is given by
$((x,\cdot), \, (y, \cdot ))$,
respectively
$((x,\cdot), (z,\cdot ))$
and
$((y,\cdot),(z,\cdot))$.
\enddemo

\demo{Proof}
We know that the covering $N$, defined by
$(x+y, \cdot)$
has 
$H_1$
as in the diagram

\midspace{3cm}
$   $\newline
So, according to 12.1, we have a map
$\pi_1(M) \rightarrow Q_{2^n}$,
where
$n \ge 3$.
The image of
$H_1(N)_{(2)}$
in
$H_1(M)_{(2)}$
is
$ker \, H_1(M) \overset {(x+y,\cdot)} \to \longrightarrow \Z/2 \Z$,
which is generated by 
$x+y$
and $z$. By discussion of $\qquad$, $u$ projects to
$x+y$,
so
$v$ and $w$ project either to $z$, or to
$x+y +z$.
Now, the induced map
$H_1(M) \rightarrow \Z/2 \Z \oplus \Z/2 \Z$
which we want to understand is just the quotient by the subspace generated by the image of $v$ or $w$. Suppose $v$ and $w$ project to
$x + y + z$,
then the map
$H_1(M) \rightarrow \Z/2\Z \oplus \Z/2 \Z$
is given by
$((x+y, \cdot ),(x+z, \cdot))$.
If
$n=3$
then by the lemma below we should have
$(x+y)^2 + (x+z)^2 + (x+y)(x+z) =0$,
however, this expression equals
$x^2 + y^2 + z^2$
by the previous Theorem. So
$n=4$,
but then we should have
$(x+y)(x+z) =0$
(or maybe
$(x+z)(y+z) =0$
or
$(x+y)(y+z) =0$).
But
$(x+y)(x+z) = x^2 \ne 0$,
so this is also impossible. Hence $u$ and $v$ project to $z$ and the map
$H_1(M) \rightarrow \Z/2\Z \oplus \Z/2 \Z$
indeed looks like
$((x,\cdot),(y,\cdot ))$.
Since
$x^2 + y^2 + xy \ne 0$,
we see that
$n \ge 4$. 
This proves the Corollary.
\enddemo

\proclaim{Lemma (15.5)}
(a) For any generators
$x,y$
of
$H^1(Q_8.\F_2)$
we have
$x^2 + y^2 + xy =0$. \newline
(b) There are generators
$x,y$
of
$H^1(Q_{2^n},\F_2), \, n \ge 4$,
such that
$xy =0$.
\endproclaim

\demo{Proof}
Consider the extension
$0 \rightarrow \Z/2\Z \rightarrow Q_8 \rightarrow \Z/2 \Z \oplus \Z/2 \Z \rightarrow 0$.
There is a 
$\Z/3\Z$
-action in this exact sequence, so the class of extensions should be
$\Z/3\Z$
-invariant. But the only invariant element in
$H^2(\Z/2 \Z \oplus \Z/2 \Z , \F_2)$
is 
$x^2+ y^2 + xy$,
as mentioned in $\qquad \qquad$, which proves (a). To prove (b) we notice that there is a map
$Q_{2^n} \twoheadrightarrow D_{2^{n-1}}$,
which induces isomorphisms in
$H^1$.
Similarly,
$D_{2^{n-1}}$
maps on
$D_8$.
Now, the extension class of
$0 \rightarrow \Z/2 \Z \rightarrow D_8 \rightarrow \Z/2\Z \oplus \Z/2 \Z \rightarrow 0$
is
$xy$
in some basis [B]. This proves (b).
\enddemo

\head{16. Homology of coverings for $H_1(M) = \Z/2 \Z \oplus \Z/4\Z$}\endhead

\demo{16.1}
In this chapter we do actual computation of the homology of all
$C_2$
-coverings in the case
$H_1(M)_{(2)} = \Z/2\Z \oplus \Z/4\Z$,
the smallest nontrivial case. This computation is based heavily on the results of the previous chapter. Here is our main result.
\enddemo

\proclaim{Theorem (16.1)}
Suppose
$H_1(M)_{(2)} = \Z/2\Z \oplus \Z/4\Z$
and $M$ is not virtually Haken. Then the 2-torsions in three
$C_2$
-coverings of $M$ are
$\Z/2 \Z \oplus \Z/2 \Z \oplus \Z/2\Z$,
$\Z/4\Z \oplus \Z/4\Z$
and
$\Z/2\Z \oplus \Z/2^{n+1}\Z, \, n\ge 1$.
\endproclaim

\demo{Proof}
By assumption the first homology of $M$ looks like

\midspace{2cm}
$   $\newline
Now, let
$U = (2u,\cdot )$
and
$V = (v,\cdot)$
be a natural basis in
$H^1(M,\F_2)$.
Then
$U^2 =0$.
Since
$UV \cdot U =0$
and
$UV \cdot V = UV^2 = U^2V =0$
we see that
$U \cdot V =0$
by Poincar\'e duality. On the other hand,
$V^2 = \beta V \ne 0$,
and therefore
$V^3 \ne 0$
again by Poincar\'e duality. So the ring
$H^*(M,\F_2)$
looks like
$U^2=0, \, UV =0, \, V^3 =1$.

Let $N \rightarrow M$
be a 
$C_2$
-covering, defined by
$(2u, \cdot) = U$.
By the Shrinking Theorem $\qquad \qquad$,
$W = H_1(N)_{(2)}$
looks like
$\Z/2\Z \oplus S$.
Moreover, $S$ should be a closed module with
$S^{C_2} \approx \Z/2\Z$.
Hence $S$ looks like

\midspace{2cm}
$   $\newline
Now if
$n > 2$,
then $S$ is a cyclic group
$\Z/2^n\Z$
with the action by multiplication by
$\pm (2^{n-1}-1)$,
which is impossible by the Theorem 12.2 (6). So
$n = 2$
and $W$ looks like

\midspace{2cm}
$   $\newline
where $x$ is anisotropic, but $p$ and $q$ may be isotropic. However, since
$\zeta p = q$
we claim that
$p+q$
is isotropic, because
$(p+q, p+q) = (p,p) + (q,q) = 2(p,p) =0$.
We now consider the covering
$Q \rightarrow N$
defined by
$(p+q)$.
The action corresponding to this covering in
$R = H_1(Q)_{(2)}$
will be called
$\eta$,
and we keep the notation
$\zeta$
for the induced action in $R$. So, as a 
$\eta$-module, $R$ looks like

\midspace{2cm}
$  $\newline
by the proof of Theorem 15.4. We wish now to understand the possibilities for
$\zeta$-action.
\enddemo

\demo{Case A}
$t$ may be chosen so that the cyclic group $T$ generated by it is
$\zeta$
-invariant. By the Theorem 12.2 (b) the action of
$\zeta$
in $T$ must be multiplication by
$(-1)$.
But then the free involution
$\zeta \eta$
will act trivially on $T$ which is impossible by the same Theorem. So the case $A$ is impossible.
\enddemo

\demo{Case B}
Since $t$ projects (to $W$) to
$p + q$,
which is
$\zeta$-invariant, we must have
$\zeta t = (2m+1)t +$
(something that projects to zero). In other words,
$\zeta t = (2l+1)t + (r+s)$.
Now, for
$\zeta r$
we have four possibilities:
\roster
\item"(i)" $\zeta r =r$ \newline
\item"(ii)" $\zeta r = r+2^{n-1}t$ \newline
\item"(iii)" $\zeta r =s$ \newline
\item"(iv)" $\zeta r = s + 2^{n-1}t$
\endroster
\enddemo

Since
$\zeta$
and
$\eta$
commute, we have
$\zeta s=s$
in case (i). But then the group generated by $r$ and $s$ is
$\zeta$
-invariant. Since the action of
$\zeta$
is orthogonal, its complement, that is $T$, is also
$\zeta$-invariant and we are brought back to the case A. This shows (i) is impossible. (iii) is reduced to $i$ by changing
$\zeta$
to
$\zeta \eta$.
Likewise (ii) is equivalent to (iv), so we will only study (iv). In this case we have
$\zeta r = s+2^{n-1}t, \, \zeta s = r+2^{n-1}t, \, \zeta t = (2l+1) t + (r+s)$.
It follows that
$Ker \, (1-\zeta)$
is generated by
$r+s$
and
$2t$
in case
$2l+1 =1$, 
or
$r + s$
and
$2^{n-1} t$
in case
$2l + 1 = -1$,
or
$r + s$
and
$t+r$
in case
$2l+1 = 2^{n-1}+1$
or
$r+s$
and
$2^{n-1}t$
in case
$2l+1 = 2^{n-1}-1$.
In all cases
$(1+\zeta )r = r+s+2^{n-1}t$
and
$\zeta t+t = (2l+2)t + (r+s)$.
If
$2l+2 \ne 0, \, 2^{n-1}$,
then
$\zeta(2t) + 2t = 2(2l+2)t$
and
$2^{n-1}t \in Im \, (1+\zeta)$,
so also
$r+s = (r+s+2^{n-1}t) -2^{n-1}t \in Im (1+\zeta)$,
so
$(2l+2)t \in Im (1+\zeta)$.
That means
$Ker \, (1-\zeta)= Im \, (1+\zeta)$,
for
$2l+1 =1$.
If
$2l+1=-1$
then
$(1+\zeta )t = r+s$,
so
$r+s \in Im \, (1+\zeta)$
and since
$\zeta r+r = r+s+2^{n-1}t$,
also
$2^{n-1}t \in Im \, (1+\zeta)$
so again
$Ker \, (1-\zeta) = Im \, (1+\zeta)$.
The final case
$2l+1 = 2^{n-1} \pm 1$
will be studied later.

So in all cases we had so far,
$H^1(\Z/2 \Z, H_1(Q)) =0$
for the
$\zeta$
-action. We can replace the chain of coverings
$Q \overset {\eta} \to \longrightarrow N \overset {\zeta} \to \longrightarrow M$
by
$Q \overset {\zeta} \to \longrightarrow G \overset {\eta} \to \longrightarrow M$
and then see that we have the situation of case $B$ of 3.1 which has been shown to be impossible in 3.5.

Now consider the possibility

$$\align
\zeta t &= (2^{n-1}-1) t + (r+s) \\
\zeta r &= s + 2^{n-1}t \\
\zeta s &= r + 2^{n-1}t
\endalign$$
Relabeling
$\zeta \eta$
by
$\zeta$
we have an action

$$\align
\zeta t &= (2^{n-1} + 1)t + (r+s) \\
\zeta r &= r+2^{n-1} t \\
\zeta s &= s + 2^{n-1} t
\endalign$$
Now, the kernel
$Ker \, (1-\zeta)$
is generated by
$r+2^{n-2}t$
and
$s+2^{n-2}t$
subject the relation 
$2((r+2^{n-2}t)+ (s+2^{n-2}t)) =0$,
so
$Ker \, (1-\zeta) \approx \Z/2\Z \oplus \Z/4\Z$.
Moreover,
$Im \, (1+\zeta)$
is generated by
$2^{n-1}t$
and
$(r+s)$,
so
$Im \, (1+\zeta) = \Z/2 \Z \oplus \Z/2 \Z$,
so
$H^1(\Z/2 \Z, H_1(Q)) = \Z/2\Z$,
which means that
$Q \rightarrow G$
is isotropic. Moreover, since
$H_1(Q)_{C_2} \approx H_1(Q)^{C_2} \approx \Z/2 \Z \oplus \Z/4\Z$,
we see that
$H_1(G)_{(2)}$
contains
$\Z/2\Z \oplus \Z/4\Z$
as an index 2 subgroup. Hence there are three possibilities for
$H_1(G)_{(2)}$:
\roster
\item"Case I".
$H_1(G)_{(2)} = \Z/2\Z \oplus \Z/2\Z \oplus \Z/4 \Z$ \newline
\item"Case II".
$H_1(G)_{(2)} = \Z/4 \Z \oplus \Z/4 \Z$ \newline
\item"Case III".
$H_1(G)_{(2)} \approx \Z/8 \Z \oplus \Z/2\Z$
\endroster
On the other hand,
$G \overset {\eta} \to \longrightarrow M$
is given by anisotropic element $v$, so
$H_1(G)_{(2)}$
under the
$\eta$
-action should be a direct sum of closed modules and
$(H_1(G)_{(2)})_{C_2}$
should be
$\Z/4 \Z$,
so we actually have just one closed module, which must look like

\midspace{2cm}
$   $\newline
It follows that only case $I$ can be realized, and
$m=3$,
so our module is

\midspace{2cm}
$   $\newline
which is just
$\Z/4 \Z \oplus \Z/4 \Z$
with permutation action. The remaining case

$$\align
\zeta t &= (2^{n-1}+1)t + (r+s) \\
\zeta r &= s+2^{n-1}t \\
\zeta s &= r + 2^{n-1}t
\endalign$$
is dealt in the same manner. In fact, we determined homology of \underbar{all}
$C_2$
-coverings of $M$, these are:
$\Z/2 \Z \oplus \Z/2 \Z \oplus \Z/2 \Z$
for the unique isotropic extension and
$\Z/4\Z \oplus \Z/4 \Z$
and
$\Z/2 \Z \oplus \Z/2^{n+1} \Z$
for two anisotropic extensions.


\enddocument

\end